\documentclass[a4, 11pt,  reqno]{amsart}
\usepackage{fullpage}

\usepackage{amsmath}
\usepackage{amsthm}
\usepackage{amsfonts}
\usepackage{amssymb}
\usepackage{amsaddr}
\usepackage{tikz-cd}
\usepackage[utf8x]{inputenc}
\usepackage{enumerate}
\usepackage[english]{babel}
\usepackage{esint}
\usepackage{hyperref}
\usepackage{bm}
\usepackage{bbm}
\usepackage[mathscr]{eucal}

\theoremstyle{plain}
\newtheorem{theorem}{Theorem}[section]
\newtheorem{lemma}[theorem]{Lemma}

\theoremstyle{definition}

\theoremstyle{remark}

\newtheorem{remark}[theorem]{Remark}

\numberwithin{equation}{section}

\newcommand{\vr}{\varrho}
\newcommand{\vu}{\textbf{\textup{u}}}
\newcommand{\vv}{\textbf{\textup{v}}}

\newcommand{\vm}{\textbf{\textup{m}}}

\newcommand{\ve}{\varepsilon}
\newcommand{\vS}{\mathbb{S}}

\newcommand{\bfphi}{\bm{\varphi}}

\newcommand{\Div}{{\rm div}_x}
\newcommand{\Grad}{\nabla_x}
\newcommand{\DerTime}{\partial_t}
\newcommand{\dx}{\textup{d}x}
\newcommand{\dt}{\textup{d}t}
\newcommand{\ds}{\textup{d}s}

\definecolor{falured}{rgb}{0.5, 0.09, 0.09}
\definecolor{pinegreen}{rgb}{0.0, 0.47, 0.44}
\definecolor{denim}{rgb}{0.08, 0.38, 0.74}

\def\softd{{\leavevmode\setbox1=\hbox{d}%
		\hbox to 1.05\wd1{d\kern-0.4ex{\char039}\hss}}}

\allowdisplaybreaks[1]
\makeatletter
\def\@settitle{\begin{center}%
  \baselineskip14\p@\relax
    \huge
  \@title
  \end{center}%
}
\makeatother

\hypersetup{colorlinks=true, linkcolor=blue, citecolor=red, urlcolor=blue}

\begin{document}
	
	\title{Global weak solutions to a compressible Navier--Stokes/Cahn--Hilliard system \\ 
    with singular entropy of mixing}

\author[Danica Basarić  \& Andrea Giorgini]{Danica Basari\'{c} \& Andrea Giorgini}

\address{Politecnico di Milano\\
Dipartimento di Matematica\\
Via E. Bonardi 9, I-20133 Milano, Italy\\
\href{mailto:danica.basaric@polimi.it}{danica.basaric@polimi.it},
\href{mailto:andrea.giorgini@polimi.it}{andrea.giorgini@polimi.it}}

\begin{abstract}
We study a Navier-Stokes/Cahn-Hilliard system modeling the evolution of a compressible binary mixture of viscous fluids undergoing phase separation. The novelty of this work is a free energy potential including the physically relevant Flory-Huggins (logarithmic) entropy, as opposed to previous studies in the literature, which only consider regular potentials with polynomial growth. 
Our main result establishes the existence of global-in-time weak solutions in three-dimensional bounded domains for arbitrarily large initial data.
The core contribution is the derivation of new estimates for the chemical potential and the Flory-Huggins entropy arising from a density-dependent Cahn-Hilliard equation under minimal assumptions: non-negative $\gamma$-integrable density with $\gamma>\frac32$. In addition, we prove that the phase variable, which represents the difference of the mass concentrations, takes value within the physical interval $(-1,1)$ almost everywhere on the set where the density is positive.
\end{abstract}

	\maketitle 
	
	\section{Introduction}

    We investigate the compressible Navier-Stokes/Cahn-Hilliard system
    \begin{align}
		\DerTime \vr + \Div(\vr \vu ) &= 0, \label{A1}\\
		\DerTime (\vr \vu ) + \Div(\vr \vu \otimes \vu )  + \Grad p&= \Div \vS  - \Div \left( \Grad c \otimes \Grad c - \frac{|\Grad c|^2}{2} \mathbb{I}\right)\!,  \label{A2} \\
		\DerTime(\vr c) + \Div (\vr c \vu) &= \Delta_x \mu, \label{A3}
	\\
    \vr \mu &= - \Delta_x c + \vr  \frac{\partial f(\vr, c)}{\partial c }. \label{A4}
    \end{align} 
    System \eqref{A1}--\eqref{A4} describes the motion of a mixture of two viscous and compressible fluids incorporating phase separation mechanisms. The state variables (unknowns)  in \eqref{A1}--\eqref{A4} are the density $\vr= \vr(t,x)$ of the mixture, the velocity $\vu=\vu(t,x)$ of the fluid mixture, the phase variable $c= c(t,x)$, representing the mass concentration difference of the two components, and the chemical potential $\mu=\mu(t,x)$. More precisely, from the continuum theory of mixtures, $c(t,x)=c_1(t,x)-c_2(t,x)$, where the mass concentrations are defined by
    \begin{equation}
    \label{c_i-def}
    c_i(t,x)= \lim_{|V|\to 0} \frac{m_i(V)}{m(V)}, \quad i=1,2 \quad \text{with}\quad m(V)=m_1(V) + m_2(V),
    \end{equation}
  where $m_i(V)$ denotes the mass at time $t$ of the fluid $i$ in the arbitrary control volume $V$ around $x$.
    
    The pressure $p=p(\vr, c)$ is related to the free energy density $f=f(\vr, c)$ through the constitutive law
	\begin{equation}
    \label{p-law}
		p(\vr, c) = \vr^2 \frac{\partial f(\vr, c)}{\partial \vr }, 
	\end{equation}
	and the viscous stress tensor $\vS = \vS( c, \Grad \vu )$ satisfies Newton's rheological law, 
	\begin{equation}
		\vS(c, \Grad \vu) =  \eta(c) \left( \Grad \vu + \Grad^{\top} \vu - \frac23 \Div \vu \, \mathbb{I}\right) + \lambda(c) \Div \vu \, \mathbb{I}.
	\end{equation}
	The functions $\eta= \eta(c)$ and $\lambda=  \lambda(c)$ are the shear and bulk viscosities, respectively. We assume that  $\eta, \, \lambda$ are continuously differentiable functions such that 
	\begin{equation} \label{viscosity coefficients}
		0 < \underline{\eta} \leq \eta(c) \leq \overline{\eta}, \quad 0 \leq \lambda(c) \leq \overline{\lambda} \quad \mbox{for all } c \in \mathbbm{R}.
	\end{equation} 
The free energy density is given by
	\begin{equation*}
		f(\vr, c) = f_{\rm e}(\vr)+ f_{\rm mix}( c), 
	\end{equation*}
where $f_{\rm e}(\vr)$ and $f_{\rm mix}(c)$ represent the potential and mixing contributions, respectively. In particular, we consider
	\begin{align*}
		f_{\rm e}(\vr)= \int_{1}^{\vr} \frac{p_{\rm e}(z)}{z^2} \ \textup{d}z, \quad &p_{\rm e}(\vr) = \vr^{\gamma}, \ \gamma>1, 
	\end{align*}
and
\begin{align*}
		f_{\rm mix}( c)=  F(c) & - \frac{\theta_0}{2} c^2,
\end{align*}
	where $F=F(c)$ denotes the Flory-Huggins (also known as Boltzmann-Gibbs) entropy of mixing 
	\begin{align}
    \label{FL-pot}
		F(c) = \frac{\theta}{2} \bigg((1+c) \ln (1+c) + (1-c) \ln (1-c)\bigg), \quad \mbox{for any } c\in[-1,1].
	\end{align}
	The parameters $\theta$ and $\theta_0$ satisfy the thermodynamical conditions $0< \theta < \theta_0$.
	
	We consider the system \eqref{A1}--\eqref{A4} in $(0,T) \times \Omega$, where the positive time $T$ can be arbitrarily chosen, and $\Omega \subset \mathbbm{R}^3$ denotes a bounded domain of class $C^2$. We close the system with the boundary conditions
	\begin{equation} \label{boundary conditions}
		\vu = \textbf{0}, \quad \Grad c \cdot \textbf{n} =0, \quad \Grad \mu \cdot \textbf{n} =0 \quad \text{on } (0,T)\times \partial \Omega,
	\end{equation}
    and the initial conditions
	\begin{equation}
    \label{initial conditions}
		\vr(0, \cdot) = \vr_0, \quad (\vr \vu)(0,\cdot) =\vm_0, \quad  c(0, \cdot) =c_0  \quad \text{in } \Omega.
	\end{equation}  
Two main properties of system \eqref{A1}--\eqref{A4}, subject to \eqref{boundary conditions}--\eqref{initial conditions}, are 
the mass conservation relations and the energy balance. More precisely, sufficient smooth solutions to \eqref{A1}--\eqref{A4} and \eqref{boundary conditions}--\eqref{initial conditions} satisfy
\begin{align}
		\int_{\Omega} \vr(t, \cdot) \ \dx = \int_{\Omega} \vr_{0} \ \dx, \quad
		\int_{\Omega} (\vr c )(t, \cdot) \ \dx = \int_{\Omega} \vr_{0} c_0  \ \dx,
	\end{align}
for all $t \geq 0$, and
\begin{equation} 
\label{energy balance}
E(\vr(t), \vu(t),c(t)) 
+ \int_{0}^{t} \int_{\Omega} \vS (c, \Grad \vu): \Grad \vu \ \dx \ds + \int_{0}^{t} \int_{\Omega} |\Grad \mu|^2 \ \dx \ds = 
E(\vr_0, \vu_0,c_0),
\end{equation}
for all $t \geq 0$, where the total energy of the system is      
\begin{equation*}
E(\vr, \vu, c) :=\int_{\Omega} \left( \frac12 \vr |\vu|^2 + \frac{\vr^{\gamma}}{\gamma-1}+ \vr F(c)-  \frac{\theta_0}{2}\vr c^2   + \frac12 |\Grad c|^2 \right) \ \dx.
\end{equation*}

\subsection{Background, literature review and aim of this work}

    Navier-Stokes/Cahn-Hilliard systems arise within the framework of Diffuse Interface methods, also referred to as the Phase Field approach, to model multi-phase flows undergoing phase separation. In this formulation, the interface between the two constituents (e.g. fluids) is described as a diffuse layer of finite thickness. The order parameter $c$ is allowed to smoothly vary between the pure phases $+1$ and $-1$ at the interface, with the interface identified as level sets of $c$. The evolution of the state variables, including the density and the velocity, are derived by combining the continuum theory of mixtures and an energy-based variational energy framework, which capture the multi-scale nature of the phenomena. Due to its ability to handle topological changes and singular interfacial behaviors, the Diffuse Interface approach has become a widely used and versatile tool in computational fluid dynamics. We refer the reader to \cite{AGGM, GvdZ} and the references therein for a broad overview.

    The first Navier-Stokes/Cahn-Hilliard model, also known as Model H, was introduced by Hohenberg and Halperin \cite{HH1977}, and later derived within the continuum mechanics setting by Gurtin, Polignone and Vi\~{n}als \cite{GPV1996}. Since then, Navier-Stokes/Cahn-Hilliard systems have been the subject of extensive theoretical and computational research, especially in the case of incompressible or quasi-incompressible mixtures. In particular, since the Model H relies on the assumption of constant density for the whole mixture, several generalizations have been developed in \cite{AGGrun, ADGK, B2002, DSS, EZAS, HMR, LT, SYW, SSBZ2017} to account for mixtures of two incompressible fluids with unmatched constant densities. We refer in particular to \cite{EZAS} for a comparative analysis of these models. 
    
    The progresses in modeling have been complemented by a comprehensive mathematical theory addressing well-posedness, regularity properties of solutions, long-time convergence to stationary states, existence of the global attractor, as well as nonlocal-to-local limits. Without claiming for completeness, we mention \cite{AbelsCMP, Abels2009, ADG, AGG, BOYER1999, B2001, FRI0, FRI, GGGP, GG, GGW, GMT, HKP}, where the Flory-Huggins potential \eqref{FL-pot} has been mostly considered. In addition, a non-homogeneous incompressible Navier-Stokes/Cahn-Hilliard model was studied by the second author and Temam in \cite{GioTem}, where the existence of global-in-time weak solutions is demonstrated under the assumption of strictly positive and bounded density.

    In contrast, considerably less is known concerning compressible Navier-Stokes/Cahn-Hilliard systems. A first compressible model was introduced by Lowengrub and Truskinovsky in \cite{LT} (see also \cite{EBS2025} for recent extensions to $N$-phase mixtures). Subsequently, a simplified variant was proposed by Abels and Feireisl in \cite{AbeFei}, which corresponds to the system \eqref{A1}--\eqref{A4} considered in the present work. In that paper, the existence of global-in-time weak solutions was established under the assumption that the free energy density of mixing takes the form
    \begin{equation}
        f_{\rm mix}(\varrho,c)= H(c)\ln (\varrho)+G(c), 
    \end{equation}
    with the following growth conditions on $H$ and $G$:
    $$
    - \underline{H}\leq H(c), H'(c)\leq \overline{H}, \quad 
    \underline{G}_1 c- \underline{G}_2 \leq G'(c) \leq \overline{G}(1+c) \quad \forall \, c \in \mathbbm{R},
    $$
    where $\underline{H}, \overline{H},\underline{G}_1, \underline{G}_2, \overline{G}$ are positive constants.
    More recently, Elbar and Poulain \cite{EP} proved the existence of global-in-time weak solutions to system \eqref{A1}--\eqref{A4}, with an additional friction term $\kappa(\vr, c) \vu $ in \eqref{A2}, under the conditions 
    $$
    \gamma >6 \quad \text{and} \quad f_{\rm mix}(c)=\frac14 (c^2-1)^2.
    $$    
    On the other hand, Kotschote and Zacher in \cite{KotZac} established local-in-time existence and uniqueness of strong solutions to the model proposed in \cite{LT}. Further contributions include the existence of global-in-time weak solutions in the case with dynamic boundary conditions  \cite{CFMMPP}, the low Mach number limit in \cite{ALN}, and the studies of the stationary problem in \cite{STAT1} and \cite{STAT2}. 
    
    Despite the above-mentioned contributions, a significant gap remains unsolved between the analysis of incompressible and compressible Navier-Stokes/Cahn-Hilliard systems. All existing results for the compressible case deal with regular mixing potentials with polynomial growth (cf. Remark 1.1 in \cite{AbeFei}), which are not suitable to enforce the physical constraint $c \in [-1,1]$ (cf. \eqref{c_i-def}). Consequently, the existence of \textit{physical} solutions to the compressible Navier-Stokes/Cahn-Hilliard model \eqref{A1}--\eqref{A4} with the physically relevant Flory-Huggins potential remains an open problem. 
    
    We conclude this literature review by mentioning some works related to  compressible Navier-Stokes/Allen-Cahn systems with  Flory-Huggins potential. In \cite{FeiPetRocSch},  Feireisl et al. proved the existence of global-in-time weak solutions in the case with singular pressure law in the spirit of Carnahan and Starling (namely, the density is bounded). A key achievement of their analysis is the existence of \textit{physical} solutions for the order parameter $c$, obtained through the use of the classical maximum principle, a tool that is not applicable in the fourth-order Cahn-Hilliard setting. Furthermore, in the Allen-Cahn case, the energy balance yields a direct control of the $L^2$-norm of the chemical potential $\mu$, whereas for the Cahn-Hilliard equation, only a control of the $W^{1,2}$-seminorm is available via the energy estimates (cf. \eqref{energy balance}). 
    Finally, we also recall the work of Kotschote \cite{Kot}, where the existence and uniqueness of local strong solutions were proven, and the study of traveling-wave phase boundaries carried out in \cite{Frei}. 

    The aim of this paper is to demonstrate the first result on the existence of global-in-time weak solutions to system \eqref{A1}--\eqref{A4}, subject to \eqref{boundary conditions}--\eqref{initial conditions}, under the thermodynamically relevant assumptions \eqref{p-law}--\eqref{FL-pot} with $\gamma>\frac32$. Our approach is based on a compactness method. Specifically, we first approximate the Flory-Huggins potential $F$ by a sequence of potentials $F_\varepsilon$ satisfying the structural assumptions of \cite{AbeFei}. This allows us to construct a sequence of approximating global weak solutions $(\vr_{\ve}, \vu_{\ve}, c_{\ve})$ to a regularized problem replacing $F$ with $F_\ve$ (cf. \eqref{WF2}--\eqref{WF5} and \eqref{energy inequality}).
    Next, the core of the proof consists in deriving uniform bounds for the approximating solutions $(\vr_{\ve}, \vu_{\ve}, c_{\ve})$. 
    We initially deduce the mass conservation relations \eqref{M1}--\eqref{M2} and the energy estimates  \eqref{E1}--\eqref{E5}, which are a direct consequence of the energy inequality for the approximating solutions. However, these bounds are insufficient to get compactness for the pressure $p(\vr_\ve)$, the chemical potential $\mu_\ve$ and the nonlinear term $F'_\ve(c_\ve)$. To overcome this issue, we develop a new set of estimates for the density-dependent Cahn-Hilliard equation under minimal conditions on the density, namely non-negative $\gamma$-integrable density with $\gamma>\frac32$. 
    Roughly speaking, we gain uniform estimates on $\mu_\ve$, $\vr_\ve F'_\ve(c_\ve)$ and $c_\ve$ by exploiting the control of $\nabla_x \mu_\ve$ available from the energy balance (see \cite{MIR} for the case of the Cahn-Hilliard equation without density). More precisely, our key step is to derive the uniform estimate 
    \begin{equation}
    \label{NE1}
        \left\| \int_{\Omega} \vr_{\ve} F_{\ve}'(c_{\ve}) (c_{\ve}- M_r)\ \dx \right\|_{L^2(0,T)}\leq C,
    \end{equation}
    where the constant $M_r$ is suitably chosen such that $\vr_\ve (c_\ve -M_r)$ has zero total mass. Then, since the product $F_{\ve}'(c_{\ve}) c_{\ve}$ is comparable to $|F_{\ve}'(c_\ve)|$ away from $c_\ve \approx 0$, we exploit the above bound to infer the uniform control 
    \begin{equation}
    \label{NE1-1}
    \left\| \int_{\Omega} \vr_{\ve} |F_{\ve}'(c_{\ve})| \ \dx \right\|_{L^2(0,T)}\leq C,
    \end{equation}
    which, in turn, allows us to obtain a uniform estimate of $\mu_\ve$ in the full $W^{1,2}(\Omega)$-norm, namely
    \begin{equation}
    \label{NE2}
    \| \mu_{\ve} \|_{L^2(0,T; W^{1,2}(\Omega))} \leq C.
    \end{equation}
    In addition, by exploiting the relation \eqref{A4}, we obtain the uniform bounds
    \begin{equation}
    \label{NE3}
    \| \sqrt{\vr_{\ve}} F'_{\ve}(c_{\ve}) \|_{L^2(0,T; L^2(\Omega))} \leq C \quad \text{and} \quad
    \| \Grad c_{\ve} \|_{L^2(0,T; L^{\frac{6\gamma}{\gamma +3}}(\Omega; \mathbbm{R}^3 ))} \leq C \quad \mbox{with} \quad  \frac{3\gamma}{\gamma +3} >1.
    \end{equation}
    The estimates \eqref{NE1}--\eqref{NE3} significantly improve the analysis in \cite{GioTem}, as they do not require strict positivity or boundedness of the density, and a broader class of initial data is admitted (cf. Remark \ref{remark_ID}).
    The latter estimate in \eqref{NE3} is then crucially used in Subsection \ref{Pressure-est} to obtain the fundamental estimate
    \begin{equation}
    \label{NE4}
    \| p(\vr_{\ve}) \|_{L^{q(\gamma)}((0,T) \times \Omega)} \leq C \quad \mbox{with} \quad q(\gamma):= \min \left\{ \frac53 -\frac{1}{\gamma}, \ \frac32  \right\}.
    \end{equation}
    Next, building on the above uniform estimates, we apply the Lions-Feireisl theory for the compressible Navier-Stokes equations to pass to the limit in the continuity equation \eqref{A1} and in the momentum equations \eqref{A2}. Regarding the density-dependent Cahn-Hilliard equation \eqref{A3}--\eqref{A4}, we recover the strong convergence for $c_\ve$ by adapting a compactness argument originally developed in \cite{AbeFei}. 
    Notably, the combination of \eqref{NE3} with the strong convergence of $c_\ve$ enables us to show that the limit functions $\vr$ and $c$ satisfy the \textit{physical} condition
    $$
    -1 < c < 1  \quad \mbox{a.e. in } \lbrace{ \varrho > 0 \rbrace} \subset (0,T)\times \Omega,
    $$
    which is essential for the consistency of the model. 
Finally, by exploiting the uniform estimate of $\sqrt{F''(c_\varepsilon)} \nabla c_\varepsilon$ in $L^2(0,T;L^2(\Omega;\mathbbm{R}^3))$, which turns out to be the only one involving $c_\varepsilon$ independent of $\varrho_\varepsilon$, we are also able to conclude that 
    $$
    -1 \leq c \leq 1  \quad \mbox{a.e. in } (0,T)\times \Omega.
    $$ 

%
%

\subsection{The main result} 
We are now ready to state the main result of the paper.
	\begin{theorem} 
    \label{Main Result}
		Let $T>0$ be chosen arbitrarily, let $\Omega \subset \mathbbm{R}^3$ be a bounded domain of class $C^2$, and let
		\begin{equation*}
			\gamma > \frac32.
		\end{equation*}
		Assume that the initial conditions $(\vr_0, \vm_0, c_0)$ satisfy 
		\begin{align} 
			&\vr_0 \geq 0 \ \ \text{a.e. in } \Omega, \quad 
            -1 \leq  c_0 \leq 1  \ \ \text{a.e. in } \Omega, \label{I1}
            \\  
            E_0:=&\int_{\Omega} \left( \frac12 \frac{|\vm_0|^2}{\vr_0}+ \vr_0 f(\vr_0, c_0) + \frac12|\Grad c_0|^2 \right) \dx < \infty, \label{I2}
		\end{align}
		and the total and relative masses 
		\begin{equation*}
			M:= \int_{\Omega} \vr_0\ \dx>0  , \quad M_c:= \int_{\Omega} \vr_0 c_0  \ \dx
		\end{equation*}
		comply with the condition
		\begin{equation} \label{constant K}
			M_r:= \frac{M_c}{M} \in (-1,1).
		\end{equation}
		Then, the Navier-Stokes/Cahn-Hilliard system \eqref{A1}--\eqref{boundary conditions} admits a weak solution $(\vr, \vu, c)$ with chemical potential $\mu$ in the following sense:
		\begin{itemize}
			\item[(i)]  \textbf{Regularity class}: the weak solution $(\vr, \vu, c, \mu)$ satisfies
			\begin{equation*}
			\begin{aligned}
					\vr \geq 0, \ \  -1\leq c\leq 1 \ &\mbox{a.e. in } (0,T) \times \Omega, \\
					 -1 < c < 1 \  &\mbox{a.e. in } \{ \vr >0 \} \subset (0,T) \times \Omega, 
			\end{aligned}
			\end{equation*}
			where 
			\begin{equation*}
				\{ \vr >0 \}:= \{ (t,x) \in (0,T) \times \Omega \ : \ \vr(t,x) >0  \},
			\end{equation*}
            and
            \begin{equation*}
				\begin{aligned}
					(\vr, \vr \vu, \vr c ) &\in C_{\rm weak}([0,T]; L^{\gamma}(\Omega) \times L^{\frac{2\gamma}{\gamma+1}}(\Omega; \mathbbm{R}^3) \times L^{\frac{6\gamma}{6+\gamma}}(\Omega)), \\
					(\vu, c) &\in L^2(0,T; W_0^{1,2}(\Omega; \mathbbm{R}^3)) \times L^{\infty}(0,T; W^{1,2}(\Omega))\cap L^2(0,T; W^{2,\frac{2\gamma}{\gamma+1}}(\Omega)), \\
					\mu &\in L^2(0,T; W^{1,2}(\Omega)),
				\end{aligned}
			\end{equation*}
            as well as
			\begin{equation*}
				\vr F'(c)= \begin{cases}
					\vr F'(c) &\mbox{if } \vr>0, \\
					0 &\mbox{if } \vr=0,
				\end{cases}
                \quad \text{with} \quad \varrho F'(c) \in L^2(0,T; L^\frac{2\gamma}{\gamma+1}(\Omega)).
			\end{equation*}
			
			\vspace{0.05cm}
			
			\item[(ii)] \textbf{Weak formulation of the renormalized continuity equation}: given the functions $b$ and $B$ such that
			\begin{equation} \label{b and B}
				b \in (L^{\infty} \cap C)([0, \infty)), \quad B(\vr):= B(1)+ \int_{1}^{\vr} \frac{b(z)}{z^2} \ \textup{d}z, 
			\end{equation}
			the integral identity 
			\begin{equation} \label{WF2fin}
				\left[\int_{\Omega}\big( \vr B(\vr) \varphi\big)(t, \cdot) \ \dx\right]_{t=0}^{t=\tau}= \int_{0}^{\tau} \int_{\Omega} \left[\vr B(\vr)\left( \partial_t \varphi + \vu \cdot \Grad \varphi \right) -  b(\vr) \Div \vu \, \varphi  \right]\dx \dt
			\end{equation}
			holds for any $\tau \in [0,T]$ and any $\varphi \in C_c^\infty ([0,T) \times \overline{\Omega})$, with $\vr (0,\cdot)= \vr_0$. \\[0.05cm]
			
			\item[(iii)]  \textbf{Weak formulation of the balance of momentum}: the integral identity 
			\begin{equation} \label{WF3fin}
				\begin{aligned}
					\left[\int_{\Omega} (\vr \vu \cdot \bm{\varphi})(t, \cdot) \ \dx\right]_{t=0}^{t=\tau} &=\int_{0}^{\tau} \int_{\Omega} \big[ \vr \vu \cdot \DerTime \bm{\varphi} +(\vr \vu \otimes \vu) : \Grad \bm{\varphi}  \big] \, \dx \dt \\
					&+ \int_{0}^{\tau} \int_{\Omega} \Big[p(\vr) \mathbb{I}- \vS(c, \Grad \vu)\Big] : \Grad \bm{\varphi} \ \dx \dt \\
					&+ \int_{0}^{\tau} \int_{\Omega} \left( (\Grad c \otimes \Grad c ) : \Grad \bm{\varphi} - \frac{|\Grad c|^2}{2} \Div \bm{\varphi}\right) \dx \dt
				\end{aligned}
			\end{equation}
			holds for any $\tau \in [0,T]$ and any $\bm{\varphi} \in C_c^\infty ([0,T) \times \Omega; \mathbbm{R}^3 )$, with $(\vr \vu)(0,\cdot)= \vm_0$.\\[0.05cm]
			
			\item[(iv)] \textbf{Weak formulation of the equation of the concentration}: the integral identity 
			\begin{equation} \label{WF4fin}
					\left[\int_{\Omega} (\vr c \varphi)(t,\cdot) \ \dx\right]_{t=0}^{t=\tau} = \int_{0}^{\tau} \int_{\Omega} \left[\vr c \,  (\partial_t \varphi + \vu \cdot \Grad \varphi  ) - \Grad \mu \cdot \Grad \varphi \right]\dx \dt 
			\end{equation}
			holds for any $\tau \in [0,T]$ and any $\varphi \in C_c^\infty([0,T)\times \overline{\Omega})$, with $(\vr c)(0,\cdot)= \vr_{0} c_0$. \\[0.05cm]
			
			\item[(v)]  \textbf{Weak formulation of the equation of the chemical potential:} the integral identity 
			\begin{equation} \label{WF5fin}
				\int_{0}^{T} \int_{\Omega} \vr \mu \varphi \ \dx \dt = \int_{0}^{T} \int_{\Omega} \left( \vr F'(c) \varphi - \theta_0 \vr c \varphi + \Grad c \cdot \Grad \varphi\right) \dx \dt 
			\end{equation}
			holds for any $\varphi \in C_c^\infty( (0,T) \times \overline{\Omega})$. 
        \\[0.05cm]
			
			\item[(vi)] \textbf{Energy inequality}: setting, for a.e. $\tau \in (0,T)$,
			\begin{equation*}
				E(\vr(\tau), \vu(\tau), c(\tau)):=\int_{\Omega} \left( \frac12 \vr |\vu|^2 + \frac{\vr^{\gamma}}{\gamma-1}+ \vr F(c)-  \frac{\theta_0}{2}\vr c^2   + \frac12 |\Grad c|^2 \right)(\tau, \cdot) \ \dx , 
			\end{equation*}
			where 
			\begin{equation*}
				\vr F(c)= \begin{cases}
					\vr F(c) &\mbox{if } \vr>0, \\
					0 &\mbox{if } \vr=0,
				\end{cases}
			\end{equation*}
			the integral inequality 
			\begin{equation} \label{energy inequality final}
				E(\vr(\tau), \vu(\tau), c(\tau)) + \int_{0}^{\tau} \int_{\Omega} \vS (c, \Grad \vu): \Grad \vu \ \dx \dt + \int_{0}^{\tau} \int_{\Omega} |\Grad \mu|^2 \ \dx \dt \leq E_0
			\end{equation}
			holds for a.e. $\tau \in (0,T)$. 
		\end{itemize}
	\end{theorem}
    
	\begin{remark}
    \label{remark_ID}
		We point out that assumption \eqref{constant K}, namely
        $$
        \frac{M_c}{M}= \frac{\displaystyle\int_{\Omega} \vr_0 c_0  \ \dx}{\displaystyle\int_{\Omega} \vr_0 \ \dx} \in (-1,1),
        $$
        is physically grounded in the context of phase separation for a two-phase mixture. 
        Indeed, when the initial concentration exhibits a diffuse interface between the pure phases, the set
		\begin{equation*}
			A:= \Big\{ x \in \Omega \ : \ -1 < c_0(x) <1  \Big\} \cap \Big\{ x \in \Omega \ : \ \vr_0(x)>0  \Big\}
		\end{equation*}
        has positive measure, i.e. $|A|>0$. 
		Then, due to hypothesis \eqref{I1}, we have 
		\begin{equation*}
			|M_c| \leq \int_{\Omega} \vr_0 |c_0| \ \dx = \int_{A} \vr_{0} |c_0| \ \dx + \int_{A^c} \vr_{0} |c_0| \ \dx < \int_{A} \vr_{0}  \ \dx + \int_{A^c} \vr_{0}  \ \dx = M. 
		\end{equation*}
        On the other hand, the case $M_c=M$ is equivalent to $$
        \int_\Omega \vr_0 (1-c_0) \ \dx=0,
        $$
        which implies that $c_0=1$ in the set $\lbrace \vr_0 >0 \rbrace$. Similarly, if $M_c=-M$, then $c_0=-1$ in the set $\lbrace \vr_0 >0 \rbrace$. In both cases, only one fluid is initially present, and thus no phase separation should occur. Furthermore, assumption \eqref{constant K} significantly broadens the class of admissible initial conditions compared to the assumption 
        $$
        \frac{\displaystyle \int_\Omega \vr_0 c_0 \ \dx}{\vr_\ast |\Omega|} \in (-1,1), \quad \text{where} \quad 0 < \vr_\ast \leq \vr_0(x) \leq \vr^\ast \ \ \text{a.e. in } \Omega,
        $$
        for some $\vr_\ast$, $\vr^\ast \in \mathbbm{R}$, as
        required in \cite{GioTem}
	\end{remark}

	\begin{remark}
        Further remarks are in order:
\begin{enumerate}
    \item Theorem \ref{Main Result} is also valid in two-dimensional domains $\Omega$. In the latter, the value of $\gamma$ and the regularity class of solutions can be improved.

    \item More general expressions for the pressure $p_{\rm e }$ can be considered, such as
		\begin{equation*}
			p_{\rm e } \in C([0, \infty)) \cap C^1((0, \infty)), \quad p_{\rm e}(0)=0, \quad \underline{p}_1 \vr^{\gamma-1} - \underline{p}_2 \leq p_{\rm e}'(\vr) \leq \overline{p}(1+\vr^{\gamma-1}), 
		\end{equation*}
		which was taken into account in \cite{AbeFei}. In addition, a term $H(c) \ln \varrho$ can also be included in $f_{\rm mix}$ as in \cite{AbeFei}, with minor modifications in our analysis. 

    \item More general singular potentials can be considered in addition to the specific form of the Flory-Huggins potential \eqref{FL-pot}. In particular, the proof of Theorem \ref{Main Result} remains valid for any function $F: [-1,1] \to \mathbbm{R}$ in the class 
    $C([-1,1])\cap C^2(-1,1)$ such that 
    $$
    \lim_{s\to \pm 1} F'(s)= \pm \infty, \quad \text{and} \quad F''(s) \geq \alpha >0, \quad \forall \, s \in (-1,1).
    $$
\end{enumerate}
     \end{remark}

\section{Functional setting}

Let $X$ be a real Banach space. We denote by $X^*$ its dual space. For $T>0$ and $p\in[1,\infty]$, the Lebesgue space $L^p(I;X)$
denotes the set of all strongly measurable functions $f: I\to X$ that are $p$
-integrable/essentially bounded. The set of continuous functions $f: [0,T]\to X$ is denoted by $C([0,T];X)$, endowed with the supremum norm. The space $C_{\rm weak}([0,T];X)$ consists of all functions $f \in L^\infty(0,T;X)$ such that the map $t\in [0,T] \mapsto \langle \phi, f(t) \rangle_{X^* \times X}$ is continuous for all $\phi \in X^*$.

Let $\Omega $ be a bounded domain in $\mathbbm{R}^{3}$. The Sobolev spaces of functions $u:\Omega \rightarrow \mathbbm{R}$ and of vector fields $\mathbf{u}:\Omega \rightarrow \mathbbm{R}^{3}$
are denoted by $W^{k,p}(\Omega )$ and $W^{k,p}(\Omega ;\mathbbm{R}^{3})$, respectively, where $k\in \mathbb{N}$ and $1\leq p\leq \infty $. The notations $\Vert \cdot \Vert
_{W^{k,p}(\Omega )}$ and $\Vert \cdot \Vert
_{W^{k,p}(\Omega;\mathbbm{R}^3)}$ represent their norms.
We recall the classical embedding results connecting the Sobolev and Lebesgue spaces. 
    \begin{lemma}
    	Let $\Omega \subset \mathbbm{R}^3$ be a bounded Lipschitz domain, and let $k \geq 1$, $1\leq p \leq \infty$ be fixed. If $kp<3$, then the continuous embedding 
    	\begin{equation} \label{Sobolev Emb}
    		W^{k,p}(\Omega) \hookrightarrow L^q(\Omega),
    	\end{equation}
    	holds for any $1\leq q \leq p^*$, where $p^*:=\frac{3p}{3-kp}$. Moreover, the embedding is compact whenever $q\neq  p^*$. 
    \end{lemma}

We will mainly use the aforementioned lemma for $k=1$ and $p=2$, i.e.
	\begin{equation} \label{SE1}
		W^{1,2}(\Omega) \hookrightarrow L^q(\Omega) \quad \mbox{for any }1\leq q\leq 6, 
	\end{equation}
	with its dual version 
	\begin{equation} \label{SE2}
		L^q(\Omega) \hookrightarrow (W^{1,2}(\Omega))^* \quad \mbox{for any } q \geq \frac65;
	\end{equation}
	notice that the embedding in \eqref{SE2} is compact if $q\neq \frac65$.

Next, we recall the definition of total mass (spatial average) for an integrable function
$$
\langle f \rangle_{\Omega}= \frac{1}{|\Omega|} \int_\Omega f \ \dx, \quad \forall \, f \in L^1(\Omega).
$$
By the classical Poincar\'{e} inequality in bounded Lipschitz domains $\Omega \subset \mathbbm{R}^3$, there exists a positive constant $C=C(\Omega)$ such that 
\begin{equation}
    \label{PI}
    \|f \|_{L^2(\Omega)} \le C \left( \|\nabla f\|_{L^2(\Omega)} + |\langle f \rangle_{\Omega}| \right)\!, \quad \forall \, f \in H^1(\Omega).
\end{equation}
We will also make use of the following variant of the Poincar\'{e} inequality from \cite[Lemma 2.1]{AbeFei}. 

    \begin{lemma} \label{Generalized Poincare}
    Let $\Omega \subset \mathbbm{R}^3$ be a bounded Lipschitz domain. Assume that $\varrho$ is a non-negative function such that
    $$
    0<M=\int_\Omega \varrho \ \dx, \quad \int_{\Omega} \varrho^\gamma \ \dx\leq R, \quad \text{with} \quad \gamma >\frac65.
    $$
    Then, there exists a constant $C=C(\gamma, M, R)$ such that
    $$
    \left\| f- \frac{1}{|\Omega|} \int_\Omega \varrho f \ \dx \right\|_{L^1(\Omega)} \leq 
    C\| \nabla f\|_{L^2(\Omega)}
    $$
    for any $f \in W^{1,2}(\Omega)$.
    \end{lemma}

Lastly, we recall the definition and properties of the Bogovskii operator (see, e.g., \cite[Theorems III.3.3 and III.3.4]{GALDI}).
	
	\begin{lemma} \label{Bogovskii operator}
		Let $\Omega \subset \mathbbm{R}^3$ be a bounded Lipschitz domain. There exists a linear operator 
		\begin{equation*}
			\mathcal{B}: L_0^q(\Omega) \rightarrow W_0^{1,q}(\Omega; \mathbbm{R}^3), \quad 1<q<\infty,
		\end{equation*}
		such that, for any $f \in L_0^q(\Omega)= \lbrace f \in L^q(\Omega): \ \int_\Omega f \ \dx=0 \rbrace$,
		\begin{align}
			\Div \mathcal{B} (f)&= f \quad \mbox{in } \Omega, \\
			\| \mathcal{B}(f)\|_{W^{1,q}(\Omega; \mathbbm{R}^3)}& \leq C  \| f  \| _{L^q(\Omega)} \label{estimate bogovskii}.
		\end{align}
		If, in addition, $f=\Div \textbf{\textup{g}}$ with $\textbf{\textup{g}}\in L^r(\Omega; \mathbbm{R}^3)$, $1<r<\infty$, and $\textbf{\textup{g}} \cdot \textbf{\textup{n}} |_{\partial \Omega}=0$, then 
		\begin{equation} \label{estimate bog ref}
			\|  \mathcal{B}(f)\|_{L^r(\Omega; \mathbbm{R}^3)} \leq C \|  \textbf{\textup{g}}\|_{L^r(\Omega; \mathbbm{R}^3)}.
		\end{equation}
	\end{lemma}

We conclude this section by fixing the notation that will be used throughout the paper. We denote with 
	\begin{itemize}
		\item $\textbf{a}=(a_1, a_2, a_3)$ a vector in $\mathbbm{R}^3 $; 
		\item $\mathbb{A}=[a_{ij}]_{i,j=1}^3$ a square matrix in $\mathbbm{R}^{3 \times 3} $; 
		\item $\mathbb{I}=[\delta_{ij}]_{i,j=1}^3$ the unit matrix; 
		\item $\textbf{a} \cdot \textbf{b}= \sum_{i=1}^3 a_i b_i$ the scalar product of two vectors $\textbf{a}=(a_1, a_2,a_3)$ and $\textbf{b}=(b_1, b_2, b_3)$; 
		\item $\mathbb{A} : \mathbb{B} = \sum_{i,j=1}^{3} a_{ij} b_{ij}$ the scalar product of two matrices $\mathbb{A}=[a_{ij}]_{i,j=1}^3$ and $\mathbb{B}=[b_{ij}]_{i,j=1}^3$; 
		\item $\textbf{a} \otimes \textbf{b} = [a_ib_j]_{i,j=1}^3 $ the tensor product of two vectors $\textbf{a}=(a_1, a_2,a_3)$ and $\textbf{b}=(b_1, b_2, b_3)$;
		\item $\mathbbm{1}_{A}$ the characteristic function of a set $A$;
		\item $\langle f \rangle_{\Omega}$ the spatial average value of $f$ over $\Omega$; 
		\item $f_{\ve} \to f \text{ in } X$ the strong convergence of a sequence $\{f_{\ve}\}_{\ve >0}$ to $f \text{ in } X$ as $\ve \to 0$;
		\item $f_{\ve} \rightharpoonup f \text{ in } X$ the weak convergence of a sequence $\{f_{\ve}\}_{\ve >0}$ to $f \text{ in } X$ as $\ve \to 0$; 
		\item $f_{\ve }\overset{*}{\rightharpoonup} f \text{ in } X$ the weak-$*$ convergence of a sequence $\{f_{\ve}\}_{\ve >0}$ to $f \text{ in } X$ as $\ve \to 0$. 
	\end{itemize}

	\section{Approximation scheme}

    In this section, we define the approximation scheme employed to prove Theorem \ref{Main Result}. More precisely, we approximate the Flory-Huggins potential $F$ by a sequence of functions $\lbrace F_\ve \rbrace$ defined in the whole $\mathbbm{R}$, which exhibit (at most) quadratic growth at infinity in terms of $c$. Next, we exploit the result proved in \cite{AbeFei}, which guarantees the existence of weak solutions to system \eqref{A1}--\eqref{A4} and \eqref{boundary conditions}--\eqref{initial conditions}, with $F$ replaced by $F_\ve$.
	
	\subsection{Regularization of the potential energy}
	
	Following \cite{GioTem}, for any $\ve >0 $ sufficiently small we consider the following approximations of the logarithmic potential $F$, 
	\begin{align*}
		F_{\ve}(s):= \begin{cases}
			\displaystyle\sum_{k=0}^{2} \frac{1}{k!} F^{(k)}(1-\ve) \big[ s- (1-\ve)\big]^k &\mbox{for } s > 1-\ve, 
            \\[8pt] 
            F(s) &\mbox{for } -1+\ve \leq s \leq 1-\ve, 
            \\[5pt]
			\displaystyle\sum_{k=0}^{2} \frac{1}{k!} F^{(k)}(-1+\ve) \big[ s- (-1+\ve)\big]^k &\mbox{for } s < -1+\ve. 
		\end{cases}
	\end{align*}
	Clearly, $F_{\ve} \in C^2(\mathbbm{R})$ is an even function for any $\ve>0$; additionally, $F_{\ve}'$ converges uniformly to $F'$ in every compact set of $(-1,1)$ as $\ve \to 0$. Moreover, denoting 
	\begin{equation*}
		G_{\ve}(s) = f_{\rm mix, \ve}(s) := F_\ve(s)-\frac{\theta_0}{2}s^2, 
	\end{equation*}
	the function 
	\begin{equation*}
		G'_{\ve}(s):= \frac{\partial f_{\rm mix, \ve}(s)}{\partial s} = F_{\ve}'(s) -\theta_0 s
	\end{equation*}
	is such that 
	\begin{equation} \label{bounds G}
		\underline{G}_{\ve,1} s - \underline{G}_{\ve,2} \leq G'_{\ve}(s) \leq \overline{G}_{\ve}(1+s) \quad \forall \, s \in \mathbbm{R},
	\end{equation}
	for some positive constants $\underline{G}_{\ve,1}, \underline{G}_{\ve,2}, \overline{G}_{\ve}$. 
	
	\subsection{Existence of approximating solutions.}
	
	For any fixed $\ve$, the existence of weak solutions for any finite energy initial data is guaranteed by \cite[Theorem 1.2]{AbeFei}; more precisely, the following result holds. 
	
	\begin{theorem} \label{Ex Weak App}
		Let $T>0$ be chosen arbitrarily large and let $\Omega \subset \mathbbm{R}^3$ be a bounded domain of class $C^2$. Moreover, let
		$\gamma>\frac32$, 
		and let the triple of functions $(\vr_0, \vm_0, c_0)$ satisfy conditions \eqref{I1}--\eqref{I2}. 
		Then, for any fixed $\ve>0$ the Navier-Stokes/Cahn-Hilliard system \eqref{A1}--\eqref{boundary conditions} with $f(\vr, c)$ replaced with
		\begin{equation*}
			f_{\ve}(\vr, c)= f_{\rm e}(\vr) + f_{\rm mix, \ve}(c)
		\end{equation*}
		admits a weak solution $(\vr_{\ve}, \vu_{\ve}, c_{\ve})$ with chemical potential $\mu_{\ve}$ in the sense that the following holds:
		\begin{itemize}
			\item[(i)]  \textbf{Regularity class}: we have 
			\begin{equation*}
					\vr_{\ve} \geq 0 \ \ \mbox{a.e. in } (0,T) \times \Omega,
			\end{equation*}
			and 
			\begin{subequations} \label{B1}
				\begin{align}
					(\vr_{\ve}, \vr_{\ve} \vu_{\ve}, \vr_{\ve} c_{\ve} ) &\in C_{\rm weak}([0,T]; L^{\gamma}(\Omega) \times L^{\frac{2\gamma}{\gamma+1}}(\Omega; \mathbbm{R}^3) \times L^{\frac{6\gamma}{6+\gamma}}(\Omega)), \label{B1.1}\\
					(\vu_{\ve}, c_{\ve}) &\in L^2(0,T; W_0^{1,2}(\Omega; \mathbbm{R}^3)) \times L^{\infty}(0,T; W^{1,2}(\Omega)), \\
					\mu_{\ve} &\in L^2(0,T; W^{1,2}(\Omega)).
				\end{align}
			\end{subequations}
			\vspace{0.05cm}
			
			\item[(ii)] \textbf{Weak formulation of the renormalized continuity equation}: given the functions $b$ and $B$ satisfying \eqref{b and B}, the integral identity 
			\begin{equation} \label{WF2}
				\begin{aligned}
					\left[\int_{\Omega}\big( \vr_{\ve} B(\vr_{\ve}) \varphi\big)(t, \cdot) \ \dx\right]_{t=0}^{t=\tau}&= \int_{0}^{\tau} \int_{\Omega} \left[\vr_{\ve} B(\vr_{\ve})\left( \partial_t \varphi + \vu_{\ve} \cdot \Grad \varphi \right) \right]\dx \dt \\
					&- \int_{0}^{\tau} \int_{\Omega} \varphi \, b(\vr_{\ve}) \, \Div \vu_{\ve}\ \dx \dt 
				\end{aligned}
			\end{equation}
			holds for any $\tau \in [0,T]$ and any $\varphi \in C_c^\infty ([0,T) \times \overline{\Omega})$, with $\vr_{\ve} (0,\cdot)= \vr_0$. \\[0.05cm]
			
			\item[(iii)]  \textbf{Weak formulation of the balance of momentum}: the integral identity 
			\begin{equation} \label{WF3}
				\begin{aligned}
					\left[\int_{\Omega} (\vr_{\ve} \vu_{\ve} \cdot \bm{\varphi})(t, \cdot) \ \dx\right]_{t=0}^{t=\tau} &=\int_{0}^{\tau} \int_{\Omega} \big[ \vr_{\ve} \vu_{\ve} \cdot \DerTime \bm{\varphi} +(\vr_{\ve} \vu_{\ve} \otimes \vu_{\ve}) : \Grad \bm{\varphi}  \big] \, \dx \dt \\
					&+ \int_{0}^{\tau} \int_{\Omega} \Big[p(\vr_{\ve}) \mathbb{I}- \vS(c_{\ve}, \Grad \vu_{\ve})\Big] : \Grad \bm{\varphi} \ \dx \dt \\
					&+ \int_{0}^{\tau} \int_{\Omega} \left( (\Grad c_{\ve} \otimes \Grad c_{\ve} ) : \Grad \bm{\varphi} - \frac{|\Grad c_{\ve}|^2}{2} \Div \bm{\varphi}\right) \dx \dt
				\end{aligned}
			\end{equation}
			holds for any $\tau \in [0,T]$ and any $\bm{\varphi} \in C_c^\infty ([0,T) \times \Omega; \mathbbm{R}^3 )$, with $(\vr_{\ve} \vu_{\ve})(0,\cdot)= \vm_0$.\\[0.05cm]
			
			\item[(iv)] \textbf{Weak formulation of the equation of the concentration}: the integral identity 
			\begin{equation} \label{WF4}
				\left[\int_{\Omega} (\vr_{\ve} c_{\ve} \varphi)(t,\cdot) \ \dx\right]_{t=0}^{t=\tau} = \int_{0}^{\tau} \int_{\Omega} \left[\vr_{\ve} c_{\ve} \,  (\partial_t \varphi + \vu_{\ve} \cdot \Grad \varphi  ) - \Grad \mu_{\ve} \cdot \Grad \varphi \right]\dx \dt 
			\end{equation}
			holds for any $\tau \in [0,T]$ and any $\varphi \in C_c^\infty([0,T)\times \overline{\Omega})$, with $(\vr_{\ve} c_{\ve})(0,\cdot)= \vr_{0} c_0$. \\[0.05cm]
			
			\item[(v)]  \textbf{Weak formulation of the equation of the chemical potential:} the integral identity 
			\begin{equation} \label{WF5}
				\int_{0}^{T} \int_{\Omega} \vr_{\ve} \mu_{\ve} \varphi \ \dx \dt = \int_{0}^{T} \int_{\Omega} \left( 
                \vr_\ve F_\ve'(c_\ve) \varphi - \theta_0 \vr_\ve c_\ve \varphi 
                + \Grad c_{\ve}  \cdot \Grad \varphi\right)  \dx \dt 
			\end{equation}
			holds for any $\varphi \in C_c^\infty((0,T)\times \overline{\Omega})$.
			\\[0.05cm]
			
			\item[(vi)] \textbf{Energy inequality}: defining, for a.e. $\tau \in (0,T)$,
			\begin{align*}
				E_{\ve}(\vr_\ve(\tau), \vu_\ve(\tau), c_\ve(\tau))&:= \int_{\Omega} \left( \frac12 \vr_{\ve} |\vu_{\ve}|^2 + \vr_{\ve} f_{\ve}(\vr_{\ve}, c_{\ve}) + \frac12 |\Grad c_{\ve}|^2\right)(\tau, \cdot) \ \dx, \\
				E_{0,\ve}&:=\int_{\Omega} \left( \frac12  \frac{|\vm_0|^2}{\vr_0} + \vr_{0}f_{\ve}(\vr_0, c_0) + \frac12|\Grad c_0|^2 \right) \dx,
			\end{align*}
			the integral inequality 
			\begin{equation} \label{energy inequality}
				E_{\ve}(\vr_\ve(\tau), \vu_\ve(\tau), c_\ve(\tau)) + \int_{0}^{\tau} \int_{\Omega} \vS (c_{\ve}, \Grad \vu_{\ve}): \Grad \vu_{\ve} \ \dx \dt + \int_{0}^{\tau} \int_{\Omega} |\Grad \mu_{\ve}|^2 \ \dx \dt \leq E_{0,\ve}
			\end{equation}
			holds for a.e. $\tau \in (0,T)$. 
		\end{itemize} 
	\end{theorem}

	\begin{remark} \label{Larger class test function}
		The validity of \eqref{WF5} can be extended to a larger class of test functions by a density argument. 
        In view of \eqref{bounds G} and  \eqref{B1}, it is enough to take $\varphi$ in \eqref{WF5} such that
		\begin{equation*}
			\varphi \in L^2(0,T;W^{1,2}(\Omega)). 
		\end{equation*}
		In particular, for any $h>0$, any $\tau \in (0,T)$ and any $\phi \in W^{1,2}(\Omega)$, we can take 
		\begin{equation*}
			\varphi(t,x)= \mathbbm{1}_{[\tau, \tau+h]}(t) \phi(x)
		\end{equation*} 
		as test function in \eqref{WF5}, obtaining the identity 
		\begin{equation*}
			\frac1h \int_{\tau}^{\tau+h} \int_{\Omega} \left( \vr_{\ve} \mu_{\ve} \phi - 
            \vr_\ve F_\ve'(c_\ve) \phi - \theta_0 \vr_\ve c_\ve \phi 
             - \Grad c_{\ve}  \cdot \Grad \phi \right) \dx \dt =0. 
		\end{equation*}
		By Lebesgue Differentation Theorem, the limit $h \to 0$ exists for a.e. $\tau \in (0,T)$ and therefore we conclude that the integral identity 
		\begin{equation} \label{WF5bis}
			\int_{\Omega} (\vr_{\ve} \mu_{\ve})(\tau, \cdot) \phi \ \dx  = \int_{\Omega} \left( \vr_\ve F_\ve'(c_\ve) - \theta_0 \vr_\ve c_\ve \right)(\tau, \cdot) \phi \ \dx + \int_{\Omega}\Grad c_{\ve}(\tau, \cdot) \cdot \Grad \phi \ \dx 
		\end{equation}
		holds for a.e. $\tau \in (0,T)$ and any $\phi \in W^{1,2}(\Omega)$. 
	\end{remark}

	\section{Proof of Theorem \ref{Main Result}} 
	
	Let us consider the family $\{ (\vr_{\ve}, \vu_{\ve}, c_{\ve}) \}_{\ve>0}$ of weak solutions with correspondent chemical potential $\mu_{\ve}$, whose existence was stated in Theorem \ref{Ex Weak App}. We are now going to prove the following result. 
	
	\begin{lemma} \label{Lemma Conv}
		Under the hypotheses of Theorem \ref{Main Result}, passing to  suitable subsequences as the case may be, the following convergences hold:
		\begin{align}
			\vr_{\ve} \to \vr &\quad \mbox{in } C_{\rm weak} ([0,T]; L^{\gamma}(\Omega)), \label{conv rho eps}\\
			\vr_{\ve} \vu_{\ve} \to \vr \vu &\quad \mbox{in } C_{\rm weak} ([0,T]; L^{\frac{2\gamma}{\gamma+1}}(\Omega; \mathbbm{R}^3)), \label{conv rho u eps} \\ 
			\vr_{\ve}c_{\ve} \to \vr c &\quad \mbox{in } C_{\rm weak} ([0,T]; L^{\frac{6\gamma}{6+\gamma}}(\Omega)), \label{conv rho c eps} \\
			c_{\ve} \overset{*}{\rightharpoonup} c &\quad \mbox{in } L^{\infty}(0,T; W^{1,2}(\Omega)), \label{conv c eps} \\
			c_{\ve} \to c &\quad \mbox{in } L^2(0,T; W^{1,2}(\Omega)), \label{strong conv c eps} \\
			\vu_{\ve} \rightharpoonup \vu &\quad \mbox{in } L^2(0,T; W^{1,2}(\Omega; \mathbbm{R}^3)), \label{conv u eps} \\
			\mu_{\ve} \rightharpoonup \mu &\quad \mbox{in } L^2(0,T; W^{1,2}(\Omega)), \label{conv mu eps}\\ 
			\vr_{\ve} \mu_{\ve} \rightharpoonup \vr \mu &\quad \mbox{in } L^2(0,T; L^{\frac{6\gamma}{\gamma +6}}(\Omega)), \label{conv rho mu eps} \\
			\vr_{\ve}c_{\ve} \vu_{\ve} \rightharpoonup \vr c \vu & \quad \mbox{in } L^2(0,T; L^{\frac{3\gamma}{\gamma+3}}(\Omega; \mathbbm{R}^3)), \label{conv rho c u eps}\\
			\vr_{\ve} \vu_{\ve} \otimes \vu_{\ve} \rightharpoonup \vr \vu \otimes \vu &\quad \mbox{in } L^2(0,T; L^{\frac{6\gamma}{4\gamma+3}}(\Omega; \mathbbm{R}^{3\times 3})), \label{conv rho u times u eps}\\
			\vS(c_{\ve}, \Grad \vu_{\ve}) \rightharpoonup \vS(c, \Grad \vu) &\quad \mbox{in } L^2((0,T) \times \Omega; \mathbbm{R}^{3\times 3}), \label{conv S eps} \\
			\vr_{\ve} F'_{\ve}(c_{\ve}) \rightharpoonup \vr F'(c) \mathbbm{1}_{\{\vr>0\}} &\quad \mbox{in } L^2(0,T; L^{\frac{2\gamma}{\gamma+1}}(\Omega)),  \label{conv rho F' eps} \\
			p(\vr_{\ve}) \rightharpoonup p(\vr) &\quad \mbox{in } L^{\frac{4}{3}-\frac{1}{2\gamma}}((0,T) \times \Omega).
            \label{conv p eps}
		\end{align}
	\end{lemma}

	In order to achieve the aforementioned convergences, we need to establish some uniform bounds independent of $\ve$; the rest of the section will be devoted to this purpose. We will make use of the notation $C$ to denote a positive constant which may depend on the parameters of the system and $T$, but is independent of $\ve$.
	
	\subsection{Mass conservation}
	
	First of all, taking $b(\vr) \equiv 0$, $B(\vr) \equiv 1$ in \eqref{WF2}, and choosing $\varphi(t,x)= \psi(t)$ with $\psi \in C^{\infty}_c(0,T)$ as test function in \eqref{WF2} and \eqref{WF4}, we can deduce that the identities
	\begin{align}
		\int_{\Omega} \vr_{\ve}(\tau, \cdot) \ \dx &= \int_{\Omega} \vr_{0} \ \dx= M, \label{M1}\\
		\int_{\Omega} (\vr_{\ve}c_{\ve})(\tau, \cdot) \ \dx &= \int_{\Omega} \vr_{0} c_0  \ \dx= M_c \label{M2}
	\end{align}
	hold for any $\tau \in [0,T]$ in light of \eqref{B1.1}. 
	
	\subsection{Energy estimates}
	
	 From hypothesis \eqref{I1} and the fact that $F_{\ve}(c) \leq F(c)$ for any $c \in [-1,1]$, we can deduce that $E_{0, \ve} \leq E_0$ for any $\ve >0$.  Therefore, condition \eqref{I2} and the energy inequality \eqref{energy inequality} lead to the following uniform bounds: 
	\begin{align}
		\| \vr_{\ve} \|_{L^{\infty}(0,T; L^{\gamma}(\Omega))} &\leq C, \label{E1}\\
		\|\sqrt{ \vr_{\ve} }\vu_{\ve} \|_{L^{\infty}(0,T; L^2(\Omega; \mathbbm{R}^3 ))} &\leq C, \label{E2}\\
		\|\Grad c_{\ve} \|_{L^{\infty}(0,T; L^2(\Omega; \mathbbm{R}^3 ))} &\leq C, \label{E3}\\
		\|\Grad \vu_{\ve} \|_{L^2(0,T; L^2(\Omega; \mathbbm{R}^{3\times 3}) )} &\leq C, \label{E4}\\
		\|\Grad \mu_{\ve} \|_{L^2(0,T; L^2(\Omega; \mathbbm{R}^{3}) )} &\leq C. \label{E5}
	\end{align}
	Notice that we have used condition \eqref{viscosity coefficients} and the Korn inequality to deduce estimate \eqref{E4}. The latter, combined with the standard Poincar\'{e} inequality, implies the uniform bound
	\begin{equation} \label{E4bis}
		\| \vu_{\ve}\|_{L^2(0,T; W^{1,2}(\Omega; \mathbbm{R}^3))} \leq C.
	\end{equation}
	Owing to \eqref{M1}, \eqref{M2}, \eqref{E1} and \eqref{E3}, an application of Lemma \ref{Generalized Poincare} yields
	\begin{equation} \label{E6}
		\| c_{\ve}\|_{L^{\infty}(0,T; W^{1,2}(\Omega))} \leq C.
	\end{equation}
	Moreover, from the estimates \eqref{E1} and \eqref{E2}, we deduce that
	\begin{equation} \label{E7}
		\|\vr_{\ve} \vu_{\ve} \|_{L^{\infty}(0,T; L^{p_1}(\Omega; \mathbbm{R}^3 ))} \leq C \quad \mbox{with} \quad  p_1:=\frac{2\gamma}{\gamma+1} >\frac65.
	\end{equation}
	Consequently, from \eqref{E4bis} and the Sobolev embedding \eqref{SE1}, we get 
	\begin{equation} \label{E8}
		\| \vr_{\ve} \vu_{\ve} \otimes \vu_{\ve} \|_{L^2(0,T; L^{p_2}(\Omega; \mathbbm{R}^{3\times 3}))} \leq C \quad \mbox{with} \quad  p_2:= \frac{6\gamma}{4\gamma+3} >1. 
	\end{equation}
	Similarly, from estimates \eqref{E1} and \eqref{E6}, we have 
	\begin{equation} \label{E9}
		\| \vr_{\ve} c_{\ve} \|_{L^{\infty}(0,T; L^{p_3}(\Omega))} \leq C \quad \mbox{with} \quad p_3:= \frac{6\gamma}{\gamma+6} > \frac65, 
	\end{equation}
	which leads, along with \eqref{E4bis}, to
	\begin{equation} \label{E10}
		\| \vr_{\ve} c_{\ve} \vu_{\ve} \|_{L^2(0,T; L^{p_4}(\Omega; \mathbbm{R}^3))} \leq C \quad \mbox{with} \quad p_4:= \frac{3\gamma}{\gamma+3} >1. 
	\end{equation}
	
	\subsection{Estimates of chemical and free energy potentials}
	\label{Novel}
    
	\begin{itemize}
		\item[(i)] First of all, from the conservation of the total \eqref{M1} and relative \eqref{M2} masses, and  from our definition \eqref{constant K} of $M_r$, the identity
		\begin{equation} \label{I3}
			\int_{\Omega} (\vr_{\ve}c_\ve)(\tau, \cdot) \ \dx - M_r \int_{\Omega}\vr_{\ve}(\tau, \cdot) \ \dx  = M_c - M_r \, M =0 
		\end{equation}
		holds for any $\tau \in [0,T]$. Therefore, since $F'_{\ve}$ is an increasing function, we have 
		\begin{equation} \label{F2}
			   \int_{\Omega} \vr_{\ve}  F_{\ve}'(c_{\ve})
			 (c_{\ve}-M_r) \ \dx = \int_{\Omega} \vr_{\ve} \left[ F_{\ve}'(c_{\ve}) -F_{\ve}'(M_r)  \right] (c_{\ve}-M_r) \ \dx \geq 0
		\end{equation}
		a.e. on $(0,T)$. Moreover, in view of Remark \ref{Larger class test function} and \eqref{E6}, for a.e. $\tau \in (0,T)$ we can choose $\phi= c_\ve(\tau, \cdot) - M_r$ as test function in \eqref{WF5bis}, obtaining that the integral identity 
		\begin{equation} \label{F1}
			\begin{aligned}
				\int_{\Omega} &\vr_{\ve} F_{\ve}'(c_{\ve})(c_{\ve}- M_r )\ \dx +\int_{\Omega} |\Grad c_\ve|^2 \ \dx \\
				&= \int_{\Omega} \vr_{\ve} \mu_{\ve} (c_\ve-M_r) \ \dx + \theta_0 \int_{\Omega} \vr_{\ve} c_{\ve}(c_\ve-M_r)\ \dx
			\end{aligned}
		\end{equation}
		holds a.e. on $(0,T)$. Using once again \eqref{I3}, we can write 
		\begin{equation*}
			\int_{\Omega} \vr_{\ve} \mu_{\ve} (c_\ve-M_r) \ \dx = \int_{\Omega} \big(\mu_{\ve}- \langle  \mu_{\ve} \rangle_\Omega \big)(\vr_{\ve}c_\ve- M_r \vr_{\ve}) \ \dx, 
		\end{equation*}
		where $\langle \mu_{\ve} \rangle_\Omega$ denotes the spatial average value of $\mu_{\ve}$. From the Sobolev embedding \eqref{SE1} and the Poincar\'e inequality \eqref{PI}, we infer that
		\begin{equation*}
			\| \mu_{\ve}- \langle  \mu_{\ve} \rangle_\Omega \|_{L^6(\Omega)}^2 
            \leq C\|  \mu_{\ve}- \langle  \mu_{\ve} \rangle_\Omega \|_{W^{1,2}(\Omega)}^2 
            \leq C \| \Grad \mu_{\ve} \|_{L^2(\Omega)}^2. 
		\end{equation*}
		Consequently, going back to \eqref{F1}, from estimates \eqref{E6}, \eqref{E9} and applying the H\"older inequality, we obtain 
		\begin{equation*}
			\begin{aligned}
				\int_{\Omega} \vr_{\ve} F_{\ve}'(c_{\ve}) (c_{\ve}- M_r)\ \dx &\leq  \| \vr_{\ve}(c_{\ve} -M_r) \|_{L^\frac65(\Omega)} \left( \| \mu_{\ve} - \langle \mu_{\ve} \rangle_\Omega \|_{L^6(\Omega )} + \theta_0 \| c_{\ve} \|_{L^6(\Omega)}\right)\\
				& \leq C \left( \| \Grad \mu_{\ve} \|_{L^2(\Omega; \mathbbm{R}^{3})} + 1 \right). 
			\end{aligned}
		\end{equation*}
		Integrating over $(0,T)$ the previous inequality squared, and exploiting \eqref{E5}, we end up with 
		\begin{equation} \label{E11}
			\left\| 	\int_{\Omega} \vr_{\ve} F_{\ve}'(c_{\ve}) (c_{\ve}- M_r)\ \dx \right\|_{L^2(0,T)}^2 \leq C \int_{0}^{T} \| \Grad \mu_{\ve} \|_{L^2(\Omega; \mathbbm{R}^{3 })}^2 \ \dt + C \leq C.
		\end{equation}
	
		\item[(ii)]  Due to condition \eqref{constant K}, we can fix $K_0>0$ such that 
        $$|M_r|< K_0<1
        $$  
        and, for a.e. $\tau \in (0,T)$, consider the sets 
		\begin{align*}
			\Omega_0(\tau) &:= \{ x \in \Omega \ | \ | c_{\ve}(\tau, x) |< K_0  \}, \\
			\Omega_1(\tau) &:= \{ x \in \Omega \ | \  c_{\ve}(\tau, x) \geq  K_0  \}, \\
			\Omega_2(\tau) &:= \{ x \in \Omega \ | \  c_{\ve}(\tau, x) \leq  -K_0 \}. 
		\end{align*}
		Notice that 
		\begin{itemize}
			\item if $x \in \Omega_1(\tau)$, 
			\begin{equation*}
				c_{\ve}(\tau, x) - M_r \geq K_0-M_r >0 \quad \Rightarrow \quad \frac{c_{\ve}(\tau, x)-M_r}{K_0-M_r} \geq 1; 
			\end{equation*}
			\item if $x \in \Omega_2(\tau)$, 
			\begin{equation*}
				c_{\ve}(\tau, x) - M_r \leq -K_0- M_r <0 \quad \Rightarrow \quad -\frac{c_{\ve}(\tau, x)-M_r}{K_0+M_r} \geq 1. 
			\end{equation*}
		\end{itemize}
		We now write, for a.e. $\tau \in (0,T)$,
		\begin{equation*}
			\int_{\Omega} \vr_{\ve} |F_{\ve}'(c_{\ve})| \ \dx = \sum_{j=0}^{2} \int_{\Omega_j(\tau)} \vr_{\ve} |F_{\ve}'(c_{\ve})| \ \dx. 
		\end{equation*}
		Noticing that $F_{\ve}'=F_{\ve}'(s)$ is an odd and strictly increasing function on $\mathbbm{R}$, such that $F_{\ve}'(s) \geq 0$ if and only if $ s \geq 0$, we have  
		\begin{align*}
			\int_{\Omega_0(\tau)}  \vr_{\ve} |F_{\ve}'(c_{\ve})| \ \dx &\leq \max_{s\in [-K_0,K_0]}  \{|F_{\ve}'(s)|\}\int_{\Omega} \vr_\ve \ \dx \leq M \, F_{\ve}'(K_0) \leq  M  \, F'(K_0),  \\[5pt]
			\int_{\Omega_1(\tau)}  \vr_{\ve} |F_{\ve}'(c_{\ve})| \ \dx &=\int_{\Omega_1(\tau)}  \vr_{\ve} F_{\ve}'(c_{\ve}) \ \dx \leq \int_{\Omega_1(\tau)}  \vr_{\ve} F_{\ve}'(c_{\ve}) \,\frac{c_{\ve}-M_r}{K_0-M_r}   \ \dx  \\
			&\leq \frac{1}{K_0-M_r} \int_{\Omega_1(\tau)}  \vr_{\ve} F_{\ve}'(c_{\ve}) (c_{\ve}-M_r)  \ \dx, \\[5pt]
			\int_{\Omega_2(\tau)}  \vr_{\ve} |F_{\ve}'(c_{\ve})| \ \dx &=\int_{\Omega_2(\tau)}  \vr_{\ve} \big(-F_{\ve}'(c_{\ve})\big) \ \dx \leq \int_{\Omega_2(\tau)}  \vr_{\ve} \big(-F_{\ve}'(c_{\ve}) \big) \left(-\frac{c_{\ve}-M_r}{K_0+M_r} \right) \dx  \\
			&\leq \frac{1}{K_0+M_r} \int_{\Omega_2(\tau)}  \vr_{\ve} F_{\ve}'(c_{\ve}) (c_{\ve}-M_r)  \ \dx. 
		\end{align*}
        Here we have used that $|F'_\ve(s)| \leq |F'(s)|$ for any $\in (-1,1)$.
		Therefore, for a.e. $\tau \in (0,T)$, recalling also \eqref{F2}, we obtain
		\begin{align*}
			\int_{\Omega} &\vr_{\ve} |F_{\ve}'(c_{\ve})| \ \dx \\
			&\leq \frac{1}{K_0-|M_r|} \int_{\Omega_1(\tau) \cup\Omega_2(\tau)}  \vr_{\ve} F_{\ve}'(c_{\ve})  (c_{\ve}-M_r) \, \dx +M  \, F'(K_0) \\
			&= \frac{1}{K_0-|M_r|} \left|\int_{\Omega_1(\tau) \cup\Omega_2(\tau)}  \vr_{\ve} F_{\ve}'(c_{\ve})  (c_{\ve}-M_r) \, \dx \right|+M  \, F'(K_0) \\
			&= \frac{1}{K_0-|M_r|} \left|\int_{\Omega}  \vr_{\ve} F_{\ve}'(c_{\ve}) (c_{\ve}-M_r) \ \dx -\int_{\Omega_0(\tau)}  \vr_{\ve} F_{\ve}'(c_{\ve}) (c_{\ve}-M_r) \ \dx \right|  +M  \, F'(K_0) \\
			&\leq \frac{1}{K_0-|M_r|} \int_{\Omega}  \vr_{\ve} F_{\ve}'(c_{\ve}) (c_{\ve}-M_r) \ \dx + \frac{2}{K_0-|M_r|} \int_{\Omega_0(\tau)}   \vr_{\ve} |F_{\ve}'(c_{\ve})| \ \dx  +M  \, F'(K_0) \\
			&\leq \frac{1}{K_0-|M_r|} \int_{\Omega}  \vr_{\ve} F_{\ve}'(c_{\ve}) (c_{\ve}-M_r) \ \dx + \left(\frac{2}{K_0-|M_r|}+1\right) M  \, F'(K_0). 
		\end{align*}
		Integrating over $(0,T)$ the previous inequality squared, from \eqref{E11}, we infer that
		\begin{equation} \label{E12}
			\left\| \int_{\Omega} \vr_{\ve} |F_{\ve}'(c_{\ve})| \ \dx \right\|_{L^2(0,T)}^2 \leq C\left\| 	\int_{\Omega} \vr_{\ve} F_{\ve}'(c_{\ve}) (c_{\ve}-  M_r)\ \dx \right\|_{L^2(0,T)}^2 + C \leq C.
		\end{equation}
		
		\item[(iii)] We take $\phi= \mathbbm{1}_{\Omega}$ as test function in \eqref{WF5bis}, obtaining that the integral identity 
		\begin{equation*}
			\int_{\Omega} \vr_{\ve}\mu_{\ve} \ \dx = \int_{\Omega} \vr_{\ve} F_{\ve}'(c_{\ve}) \ \dx -\theta_0 \int_{\Omega} \vr_{\ve} c_{\ve} \ \dx
		\end{equation*}
		holds a.e. in $(0,T)$. From estimate \eqref{E9}, we get
		\begin{equation*}
			\left| \int_{\Omega} \vr_{\ve}\mu_{\ve} \ \dx \right| \leq \int_{\Omega} \vr_{\ve} |F_{\ve}'(c_{\ve})| \ \dx + C. 
		\end{equation*}
		Integrating over $(0,T)$ the previous inequality squared, we deduce from \eqref{E12} that
		\begin{equation*}
			\left\| \int_{\Omega} \vr_{\ve} \mu_{\ve} \ \dx \right\|_{L^2(0,T)}^2 \leq \left\| \int_{\Omega} \vr_{\ve} |F_{\ve}'(c_{\ve}))| \ \dx \right\|_{L^2(0,T)}^2 + C \leq C. 
		\end{equation*} 
		As consequence of Lemma \ref{Generalized Poincare}, we conclude that 
		\begin{equation} \label{E13}
			\| \mu_{\ve} \|_{L^2(0,T; W^{1,2}(\Omega))} \leq C.  
		\end{equation}
	
		\item[(iv)] In view of Remark \ref{Larger class test function}, we can take $\varphi= F_{\ve}'(c_{\ve}) \in L^2(0,T; W^{1,2}(\Omega))$ as test function in \eqref{WF5}, obtaining the identity
		\begin{equation*}
			\quad \int_{0}^{T} \int_{\Omega} \vr_{\ve} |F_{\ve}'(c_{\ve})|^2 \dx \dt + \int_{0}^{T} \int_{\Omega} F_{\ve}''(c_{\ve}) |\Grad c_{\ve}|^2 \ \dx \dt = \int_{0}^{T} \int_{\Omega} \vr_{\ve} F_{\ve}'(c_{\ve}) \big( \mu_{\ve}+\theta_0 c_{\ve}\big) \dx \dt. 
		\end{equation*}
		From the convexity of $F_{\ve}$, $F_{\ve}''(s)\geq 0$ for any $s \in \mathbbm{R}$. By the Young inequality, we get
		\begin{equation}
			\begin{aligned}
\int_{0}^{T} &\int_{\Omega} \vr_{\ve} |F_{\ve}'(c_{\ve})|^2 \dx \dt 
+ \int_{0}^{T} \int_{\Omega} F_{\ve}''(c_{\ve}) |\Grad c_{\ve}|^2 \ \dx \dt  \\ 
				&\leq \frac12 \int_{0}^{T} \int_{\Omega} \vr_{\ve} |F_{\ve}'(c_{\ve})|^2 \dx \dt + 
                C \int_{0}^{T} \int_{\Omega} \vr_{\ve} \mu_{\ve}^2 \ \dx \dt 
                + C \int_{0}^{T} \int_{\Omega} \vr_{\ve} c_{\ve}^2 \ \dx \dt.
                \label{rhoF'-est}
			\end{aligned}
		\end{equation}
        Exploiting the uniform estimates \eqref{E1}, \eqref{E13} and the Sobolev embedding \eqref{SE1}, it is easily seen that
        \begin{equation}
            \int_{0}^{T} \int_{\Omega} \vr_{\ve} \mu_{\ve}^2 \ \dx \dt
            \leq \int_{0}^{T} \| \vr_\ve\|_{L^\frac32(\Omega)}
            \| \mu_\ve\|_{L^6(\Omega)}^2 \ \dt
            \leq C \| \mu_{\ve} \|_{L^2(0,T; W^{1,2}(\Omega))}^2
            \leq C.
        \end{equation}
        Owing to \eqref{E6}, a similar argument holds for the last term on the right-hand side of \eqref{rhoF'-est}.
        Hence, it is straightforward to conclude that 
		\begin{equation} \label{E14}
			\| \sqrt{\vr_{\ve}} F'_{\ve}(c_{\ve}) \|_{L^2(0,T; L^2(\Omega))} \leq C
		\end{equation}	
	and
	\begin{equation} \label{E14-bis}
			\| \sqrt{F''_{\ve}(c_\ve)} \Grad c_{\ve} \|_{L^2(0,T; L^2(\Omega; \mathbbm{R}^3))} \leq C.
		\end{equation}
	
		\item[(v)] In light of \eqref{E1}, \eqref{E6}, \eqref{E13}, \eqref{E14}, we deduce that
        \begin{equation*}
        \left\| \varrho_\ve \mu_\ve - \varrho_\ve F_\ve'(c_\ve) 
        + \theta_0 \varrho_\ve c_\ve^2 \right\|_{L^2(0,T; L^{p_1}(\Omega))}\leq C, \quad \mbox{with} \quad  p_1:=\frac{2\gamma}{\gamma+1} >\frac65.
        \end{equation*}
        Hence, from the integral identity \eqref{WF5} and the elliptic regularity theory for the Laplace equation with homogeneous Neumann boundary conditions, we find that
		\begin{equation*}
			\| c_{\ve} \|_{L^2(0,T; W^{2,p_1}(\Omega))} \leq C. 
		\end{equation*} 
		Then, by the Sobolev embedding \eqref{Sobolev Emb} with $k=1$, $p=p_1$ and $q=p_1^*$ (notice that $p_1<3$), we obtain  
		\begin{equation} \label{E15}
			\| \Grad c_{\ve} \|_{L^2(0,T; L^{2p_4}(\Omega; \mathbbm{R}^3 ))} \leq C \quad \mbox{with} \quad p_4= \frac{3\gamma}{\gamma +3} >1.
		\end{equation}
		Since, by interpolation,
		\begin{equation*}
			L^{\infty}(0,T; L^2(\Omega)) \cap L^2(0,T; L^{2p_4}(\Omega)) \hookrightarrow L^{p}((0,T)\times \Omega) \quad \mbox{with} \quad p:= \frac{2(5 \gamma-3)}{3 \gamma}>2, 
		\end{equation*}
		we infer from \eqref{E6} and \eqref{E15} that
		\begin{equation} \label{E16}
			\| \Grad c_{\ve} \|_{L^p((0,T) \times \Omega; \mathbbm{R}^3 )} \leq C \quad \mbox{with} \quad p>2.
		\end{equation}
	\end{itemize}

	\subsection{Pressure estimate} 
    \label{Pressure-est}
	
	 From the energy inequality \eqref{energy inequality}, we can only deduce that 
	\begin{equation} \label{E17}
		\| p(\vr_{\ve}) \|_{L^{\infty}(0,T; L^1(\Omega))} \leq C.
	\end{equation}
	Hence, we need some supplementary pressure estimates, which can be achieved by choosing a suitable test function in the weak formulation of the balance of momentum \eqref{WF3} by following the argument devised in \cite{FP2000}.
	
	Let us consider the functions  
	\begin{equation*}
		\bfphi_{\ve}(t,x):= \psi(t) \bm{\phi}_{\ve} (t,x), \quad \psi \in C_c^{\infty}(0,T), \quad \bm{\phi}_{\ve}:= \mathcal{B} \left[ b(\vr_{\ve})- \langle b(\vr_{\ve} )\rangle_{\Omega} \right],
	\end{equation*}
	where  $\mathcal{B}$ denotes the Bogovskii operator introduced in Lemma \ref{Bogovskii operator} and $b$ is a smooth bounded function of the type $b(z)=z^{\nu}$, where the exponent $\nu>0$ will be fixed later on. To be more precise, $\bfphi_{\ve}$ is an admissible test function upon replacing $\bm{\phi}_{\ve}$ by $\bm{\phi}_{\ve} * \xi_{\alpha}$, where $\{\xi_{\alpha}\}_{\alpha>0}$ is a suitable family of mollifiers in the time variable. As this procedure is fairly standard, we omit the details below.
	We now take $\bfphi_{\ve}$ as test function in the weak formulation of the balance of momentum \eqref{WF3}, to get the identity 
	\begin{equation} \label{P1}
		\int_{0}^{T} \psi  \int_{\Omega} p(\vr_{\ve}) \, b(\vr_{\ve}) \ \dx \dt = \sum_{i=1}^{7} I_{\ve,i},
	\end{equation}
	where
	\begin{align*}
		I_{\ve,1} &:=  \frac{1}{|\Omega|} \int_{0}^{T} \psi  \left(\int_{\Omega} p(\vr_{\ve}) \ \dx\right)\left( \int_{\Omega} b(\vr_{\ve}) \ \dx\right) \dt, \\
		I_{\ve,2} &:= \int_{0}^{T} \psi  \int_{\Omega} \vr_{\ve} \vu_{\ve} \cdot  \mathcal{B} \big[ \Div \big( b(\vr_{\ve}) \vu_{\ve}\big) \big] \ \dx \dt,  \\
		I_{\ve,3} &:= \int_{0}^{T} \psi  \int_{\Omega} \vr_{\ve} \vu_{\ve} \cdot  \mathcal{B} \left[ \big(\vr_{\ve} b'(\vr_{\ve})-b(\vr_{\ve})\big) \Div \vu_{\ve} - \left\langle \big(\vr_{\ve} b'(\vr_{\ve})-b(\vr_{\ve})\big) \Div \vu_{\ve} \right\rangle_\Omega \right]  \dx \dt,  \\
		I_{\ve,4} &:= -\int_{0}^{T} \psi  \int_{\Omega} (\vr_{\ve} \vu_{\ve} \otimes \vu_{\ve}) :  \Grad \mathcal{B} \left[ b(\vr_{\ve})- \langle b(\vr_{\ve} )\rangle_\Omega \right]  \dx \dt , \\
		I_{\ve,5} &:=  \int_{0}^{T} \psi  \int_{\Omega} \vS(c_{\ve}, \Grad \vu_{\ve}) :  \Grad \mathcal{B} \left[ b(\vr_{\ve})- \langle b(\vr_{\ve} )\rangle_\Omega \right] \dx \dt, \\
		I_{\ve,6} &:= - \int_{0}^{T} \psi  \int_{\Omega} \left(\Grad c_{\ve} \otimes \Grad c_{\ve} -\frac12 |\Grad c_{\ve}|^2 \mathbb{I} \right)  :  \Grad \mathcal{B} \left[ b(\vr_{\ve})- \langle b(\vr_{\ve} )\rangle_\Omega \right] \dx \dt, \\
		I_{\ve,7} &:=- \int_{0}^{T} \partial_t \psi  \int_{\Omega} \vr_{\ve} \vu_{\ve} \cdot  \mathcal{B} \left[ b(\vr_{\ve})- \langle b(\vr_{\ve} )\rangle_\Omega \right]  \dx \dt.
	\end{align*}
	We proceed by estimating each integral $I_{\ve,i}$, $i=1,\dots,7$.
	\begin{itemize}
		\item[(i)] From \eqref{E17}, we easily find that
		\begin{equation*}
			|I_{\ve,1}| \leq C \| \psi\|_{L^{\infty}(0,T)} \| p(\vr_{\ve}) \|_{L^{\infty}(0,T; L^1(\Omega))} \| \vr_{\ve}\|_{L^\nu((0,T)\times \Omega)}^{\nu}.
		\end{equation*}
		\item[(ii)] From \eqref{estimate bog ref} and estimates \eqref{E1}, \eqref{E4bis}, we have 
		 \begin{align*}
			|I_{\ve,2}| &\leq C \| \psi\|_{L^{\infty}(0,T)} \| \vr_{\ve} \vu_{\ve} \|_{L^2(0,T; L^{p_3}(\Omega; \mathbbm{R}^3))} \left\| \mathcal{B} \big[ \Div \big( b(\vr_{\ve}) \vu_{\ve}\big) \big] \right\|_{L^2(0,T; L^{\frac{p_3}{p_3-1}}(\Omega; \mathbbm{R}^3))} \\
			&\leq C  \left\|   b(\vr_{\ve}) \vu_{\ve} \right\|_{L^2(0,T; L^{\frac{6 \gamma}{5 \gamma-6}}(\Omega; \mathbbm{R}^3))} \\
			&\leq C \left\|  \vu_{\ve} \right\|_{L^2(0,T; L^6(\Omega; \mathbbm{R}^3))} \| \vr_{\ve}\|_{L^{\infty}(0,T; L^{\frac{3\gamma  \nu}{2\gamma-3}}(\Omega))}^{\nu}.
		\end{align*}
		\item[(iii)] Define for simplicity $f_{\ve}= \big(\vr_{\ve} b'(\vr_{\ve})-b(\vr_{\ve})\big) \Div \vu_{\ve}$.
        We infer from the Sobolev embedding \eqref{Sobolev Emb} with $k=1$ and $p=\frac{6 \gamma}{5 \gamma-3}$, \eqref{estimate bogovskii}, and estimates \eqref{E1}, \eqref{E4bis}, \eqref{E7} that
		\begin{align*}
			|I_{\ve,3}| & \leq C \| \psi\|_{L^{\infty}(0,T)} 
            \| \vr_{\ve} \vu_{\ve} \|_{L^\infty(0,T; L^{p_1}(\Omega; \mathbbm{R}^3))}
            \left\| \mathcal{B} \big[ f_{\ve} - \langle f_{\ve} \rangle_\Omega \big] \right\|_{L^1(0,T; L^{\frac{p_1}{p_1-1}}(\Omega; \mathbbm{R}^3))} \\
			& \leq C  \left\| \mathcal{B} \big[ f_{\ve} - \langle f_{\ve} \rangle_\Omega \big] \right\|_{L^1(0,T; W_0^{1,\frac{6\gamma}{5\gamma-3}}(\Omega; \mathbbm{R}^3))} \\
			& \leq C \left\| \big(\vr_{\ve}b'(\vr_{\ve})-b(\vr_{\ve})\big) \Div \vu_{\ve} \right\|_{L^1(0,T; L^{\frac{6\gamma}{5\gamma-3}}(\Omega))} \\
			&\leq C \| \Div \vu_{\ve} \|_{L^2((0,T)\times \Omega)} \| \vr_{\ve}\|_{L^{2\nu}(0,T; L^{\frac{6\gamma \nu}{2\gamma-3}}(\Omega))}^{\nu}.
		\end{align*}
         
		\item[(iv)] From \eqref{estimate bogovskii} and estimates \eqref{E1}, \eqref{E4bis}, we obtain
	    \begin{align*}
			|I_{\ve,4}|&\leq C  \| \psi\|_{L^{\infty}(0,T)} \| \vr_{\ve} \vu_{\ve} \otimes \vu_{\ve}\|_{L^1(0,T; L^{p_4}(\Omega; \mathbbm{R}^{3\times 3}))} \\
			&\hspace{0.5cm} \times \| \Grad \mathcal{B} \left[ b(\vr_{\ve})- \langle b(\vr_{\ve} )\rangle_\Omega \right] \|_{L^{\infty}(0,T; L^{\frac{p_4}{p_4-1}}(\Omega; \mathbbm{R}^{3\times 3}))} \leq  \| \vr_{\ve} \|_{L^{\infty}(0,T; L^{\frac{3\gamma \nu}{2\gamma-3}}(\Omega))}^{\nu}.
		\end{align*}
		\item[(v)] From \eqref{estimate bogovskii}, hypothesis \eqref{viscosity coefficients} and estimate \eqref{E4bis}, we have 
		\begin{align*}
			|I_{\ve, 5}|&\leq C \| \psi\|_{L^{\infty}(0,T)} \| \vS (c_{\ve}, \Grad \vu_{\ve})\|_{L^2((0,T) \times \Omega; \mathbbm{R}^{3\times 3})} \\
			&\hspace{0.5cm} \times \| \Grad \mathcal{B} \left[ b(\vr_{\ve})- \langle b(\vr_{\ve} )\rangle_\Omega \right] \|_{L^2((0,T) \times \Omega; \mathbbm{R}^{3\times 3})}  \leq C \| \vr_{\ve} \|_{L^{2\nu}((0,T) \times \Omega)}^{\nu}.
		\end{align*}
		\item[(vi)] From \eqref{estimate bogovskii} and estimate \eqref{E15}, we get
		\begin{align*}
			|I_{\ve, 6}|&\leq C \| \psi\|_{L^{\infty}(0,T)} \left\| |\Grad c_{\ve}|^2 \right\|_{L^1(0,T; L^{p_4}(\Omega; \mathbbm{R}^{3\times 3}))}\\
			&\hspace{0.5cm} \times \| \Grad \mathcal{B} \left[ b(\vr_{\ve})- \langle b(\vr_{\ve} )\rangle_\Omega \right] \|_{L^{\infty}(0,T; L^{\frac{p_4}{p_4-1}}(\Omega; \mathbbm{R}^{3\times 3}))} \leq C \| \vr_{\ve} \|_{L^{\infty}(0,T; L^{\frac{3\gamma \nu}{2\gamma -3}}(\Omega))}.
		\end{align*}
		\item[(vii)] From the Sobolev embedding \eqref{Sobolev Emb} with $k=1$ and $p=\frac{6 \gamma}{5 \gamma-3}$, \eqref{estimate bogovskii} and estimate \eqref{E7}, we deduce that
		\begin{align*}
			|I_{\ve,7}| &\leq C \| \partial_t \psi \|_{L^1(0,T)} \| \vr_{\ve} \vu_{\ve} \|_{L^{\infty}(0,T; L^{p_1}(\Omega; \mathbbm{R}^3))} \left\|  \mathcal{B}_{\ve} \left[ b(\vr_{\ve})- \langle b(\vr_{\ve} )\rangle_\Omega \right] \right\|_{L^{\infty}(0,T; L^{\frac{p_1}{p_1-1}}(\Omega; \mathbbm{R}^3))} \\
			&\leq C \| \partial_t \psi \|_{L^1(0,T)}\left\|  \mathcal{B}_{\ve} \left[ b(\vr_{\ve})- \langle b(\vr_{\ve} )\rangle_\Omega \right] \right\|_{L^{\infty}(0,T; W_0^{1,\frac{6\gamma}{5\gamma-3}}(\Omega; \mathbbm{R}^3))} \\
			&\leq C \| \partial_t \psi \|_{L^1(0,T)}\| \vr_{\ve}\|_{L^{\infty}(0,T; L^{\frac{6\gamma \nu}{5\gamma-3}}(\Omega))}^{\nu}.
		\end{align*}
	\end{itemize}
	Now, to make all the above norms uniformly bounded, we have to choose 
     \begin{equation*}
		\nu \leq \min \left\{ \gamma, \ \frac23 \gamma -1, 
        \ \frac{\gamma}{3}-\frac12,
        \ \frac{\gamma}{2}, \ \frac56 \gamma -\frac12  \right\} \quad \Rightarrow \quad \nu(\gamma):= \ \frac{\gamma}{3}-\frac12.
	\end{equation*}
	If we let $\psi \rightarrow 1$ in \eqref{P1}, we can finally get
	\begin{equation*}
		\int_{0}^{T}  \int_{\Omega} p(\vr_{\ve}) \vr_{\ve}^{\nu(\gamma)} \  \dx \dt \leq C .
	\end{equation*}
	 In particular, we obtain that 
	\begin{equation} \label{E19}
		\| \vr_{\ve}\|_{L^{p(\gamma)}((0,T)\times \Omega)} \leq C \quad \mbox{with} \quad p(\gamma):= \frac43 \gamma -\frac12, 
	\end{equation}
	and therefore 
	\begin{equation} \label{E18}
		\| p(\vr_{\ve}) \|_{L^{q(\gamma)}((0,T) \times \Omega)} \leq C \quad \mbox{with} \quad q(\gamma):= \frac43 - \frac{1}{2\gamma}. 
	\end{equation}

    \begin{remark}
    If we consider the case $\frac{3}{2} < \gamma \leq 6$, then the uniform estimates for the approximating densities \eqref{E19}, and thus for the approximating pressures \eqref{E18}, can be improved. Indeed, by exploiting the Sobolev embedding \eqref{Sobolev Emb} with $k=1$ and $p=\frac{6 \gamma}{7 \gamma-6}$, \eqref{estimate bogovskii}, and estimates \eqref{E1}, \eqref{E4bis}, we have 
		\begin{align*}
			|I_{\ve,3}| & \leq C \| \psi\|_{L^{\infty}(0,T)} \| \vr_{\ve} \vu_{\ve} \|_{L^2(0,T; L^{p_3}(\Omega; \mathbbm{R}^3))} \left\| \mathcal{B} \big[ f_{\ve} - \langle f_{\ve} \rangle_\Omega \big] \right\|_{L^2(0,T; L^{\frac{p_3}{p_3-1}}(\Omega; \mathbbm{R}^3))} \\
			& \leq C  \left\| \mathcal{B} \big[ f_{\ve} - \langle f_{\ve} \rangle_\Omega \big] \right\|_{L^2(0,T; W_0^{1,\frac{6\gamma}{7\gamma-6}}(\Omega; \mathbbm{R}^3))} \\
			& \leq C \left\| \big(\vr_{\ve}b'(\vr_{\ve})-b(\vr_{\ve})\big) \Div \vu_{\ve} \right\|_{L^2(0,T; L^{\frac{6\gamma}{7\gamma-6}}(\Omega))} \\
			&\leq C \| \Div \vu_{\ve} \|_{L^2((0,T)\times \Omega)} \| \vr_{\ve}\|_{L^{\infty}(0,T; L^{\frac{3\gamma \nu}{2\gamma-3}}(\Omega))}^{\nu}.
		\end{align*}
        Therefore, we infer that 
\begin{equation*}
		\nu \leq \min \left\{ \gamma, \ \frac23 \gamma -1, 
        \ \frac{\gamma}{2}, \ \frac56 \gamma -\frac12  \right\} \quad \Rightarrow \quad \nu(\gamma):= \frac{2}{3}\gamma -1.
	\end{equation*}
    In summary, we conclude that
	\begin{equation} \label{bound rho}
		\| \vr_{\ve}\|_{L^{p(\gamma)}((0,T)\times \Omega)} \leq C \quad \mbox{with} \quad p(\gamma):= 
        \begin{cases}
			\frac53 \gamma -1 &\mbox{if } \frac32 < \gamma \leq 6, \\
			\frac{4}{3}\gamma-\frac12 &\mbox{if } \gamma >6,
		\end{cases}
	\end{equation}
	and 
	\begin{equation*}
		\| p(\vr_{\ve}) \|_{L^{q(\gamma)}((0,T) \times \Omega)} \leq C \quad \mbox{with} \quad q(\gamma):= \begin{cases}
			\frac53 -\frac{1}{\gamma} &\mbox{if } \frac32 < \gamma \leq 6, \\
			\frac{4}{3}-\frac{1}{2\gamma} &\mbox{if } \gamma >6.
		\end{cases}
        \end{equation*}
        \end{remark}
	
	\subsection{Weak convergences}
	
	First of all, convergences \eqref{conv c eps}, \eqref{conv u eps} and \eqref{conv mu eps} are a direct consequence of \eqref{E6}, \eqref{E4bis} and \eqref{E13}, respectively. From estimates \eqref{E1}, \eqref{E7} and \eqref{E9}, passing to suitable subsequences as the case may be, we deduce that
	\begin{align}
		\vr_{\ve} \overset{*}{\rightharpoonup} \vr &\quad \mbox{in } L^{\infty}(0,T; L^{\gamma}(\Omega)), \label{conv rho eps bis}\\
		\vr_{\ve} \vu_{\ve} \overset{*}{\rightharpoonup} \overline{\vr \vu} &\quad \mbox{in } L^{\infty}(0,T; L^{p_1}(\Omega; \mathbbm{R}^3)), \label{conv rho u eps bis}\\
		\vr_{\ve} c_{\ve} \overset{*}{\rightharpoonup} \overline{\vr c} &\quad \mbox{in } L^{\infty}(0,T; L^{p_3}(\Omega)), \label{conv rho c eps bis}
	\end{align}
    where the notation $\overline{f}$ stands for a weak limit in $L^1((0,T)\times \Omega)$.
	Moreover, from identity \eqref{WF2} with $b(\vr) \equiv 0$ and $B(\vr)\equiv 1$, and from estimate \eqref{E7}, we simply have  
	\begin{equation*}
		\{ \DerTime \vr_{\ve} \}_{\ve >0} \mbox{ is uniformly bounded in } L^{\infty}(0,T; [W^{1, p_1'}(\Omega)]^*),
	\end{equation*}
    where $p_1'$ is the conjugate exponent of $p_1$.
	As consequence of the Aubin-Lions lemma, we  infer that
	\begin{equation*}
		\vr_{\ve} \to \vr \quad \mbox{in } C([0,T]; [W^{\alpha, p_1'}(\Omega)]^*) \quad \mbox{for any } 0< \alpha <1, 
	\end{equation*}
	where $W^{\alpha, p}(\Omega)$ with $0<\alpha<1 $ denotes the Sobolev-Slobodeckii space. In particular, we can deduce that convergence \eqref{conv rho eps bis} can be strengthened to \eqref{conv rho eps}.
	The latter, combined with the compact Sobolev embedding \eqref{SE2}, leads to 
	\begin{equation} \label{C1}
		\vr_{\ve} \to \vr \quad \mbox{in } C([0,T]; [W^{1,2}(\Omega)]^*).  
	\end{equation}
	Consequently, from convergences \eqref{conv c eps} and \eqref{conv u eps}, we deduce that 
	\begin{equation*}
		\overline{\vr \vu} = \vr \vu, \quad \overline{\vr c} = \vr c \quad \mbox{a.e. on } (0,T) \times \Omega. 
	\end{equation*}
	Similarly, from identities \eqref{WF3}, \eqref{WF4} and the previously established estimates  (\eqref{E4bis}, \eqref{E8}, \eqref{E10}, \eqref{E13}, \eqref{E16}, \eqref{E18}), convergences \eqref{conv rho u eps bis} and \eqref{conv rho c eps bis} can be strengthened to \eqref{conv rho u eps} and \eqref{conv rho c eps}, respectively. Exploiting once again the compact Sobolev embedding \eqref{SE2}, we get from \eqref{conv rho u eps}, \eqref{conv rho c eps} and \eqref{conv u eps}, convergences \eqref{conv rho c u eps} and \eqref{conv rho u times u eps}. In addition, from \eqref{conv mu eps} and \eqref{C1}, we derive \eqref{conv rho mu eps}. 

	\subsection{Strong convergence of the concentrations}
	
	The next goal is to show that 
	\begin{equation} \label{conv grad c eps}
		\| \Grad c_{\ve} \|_{L^2((0,T) \times \Omega; \mathbbm{R}^3)} \to \| \Grad c \|_{L^2((0,T) \times \Omega; \mathbbm{R}^3)}, 
	\end{equation}
	which, along with \eqref{conv c eps}, will imply the strong convergence \eqref{strong conv c eps} of the concentrations. 
	
	The idea of the proof is essentially the same as in \cite[Section 2.6]{AbeFei} but with some simplifications; we report it for the sake of completeness.  
	\begin{itemize}
		\item From identity \eqref{WF5}, we have that
		\begin{equation} \label{CW1}
			\int_{0}^{T} \int_{\Omega} \Grad c_{\ve} \cdot \Grad \varphi \ \dx \dt = \int_{0}^{T} \int_{\Omega} \sqrt{\vr_{\ve}} f_{\ve} \, \varphi \ \dx \dt 
		\end{equation}
		holds for any $\varphi \in C^{\infty}_c ((0,T) \times \overline{\Omega})$, where, due to convergences \eqref{conv rho eps}, \eqref{conv c eps}, \eqref{conv mu eps} and estimate \eqref{E14}, 
		\begin{equation} \label{conv f eps}
			f_{\ve}:= \sqrt{\vr_{\ve}} (\mu_{\ve} -F_{\ve}'(c_{\ve}) + \theta_0 c_{\ve}) \rightharpoonup f \quad \mbox{in } L^2((0,T) \times \Omega). 
		\end{equation}
		Thus, passing to the limit in \eqref{CW1}, we deduce that the identity 
		\begin{equation} \label{CW2}
			\int_{0}^{T} \int_{\Omega} \Grad c \cdot \Grad \varphi \ \dx \dt = \int_{0}^{T} \int_{\Omega} \overline{\sqrt{\vr} f} \, \varphi \ \dx \dt 
		\end{equation}
		holds for any $\varphi \in C^{\infty}_c ((0,T) \times \overline{\Omega})$, where
		\begin{equation*}
			\sqrt{\vr_{\ve}} f_{\ve} \rightharpoonup \overline{\sqrt{\vr} f} \quad \mbox{in } L^2(0,T; L^{p_1}(\Omega)) \quad \mbox{with} \quad p_1= \frac{2\gamma}{\gamma+1} >\frac65. 
		\end{equation*}
		\item By a density argument, the class of test functions in \eqref{CW1} and \eqref{CW2} can be enlarged 
        to $L^2(0,T; W^{1,2}(\Omega))$. Therefore,  taking $\varphi=c_{\ve}$ and $\varphi= c$ in  \eqref{CW1} and \eqref{CW2}, respectively, \eqref{conv grad c eps} will follow as long as we manage to prove that 
		\begin{equation} \label{I4}
			\overline{\sqrt{\vr} f c} = \overline{\sqrt{\vr} f} \, c \quad \mbox{a.e. in } (0,T) \times \Omega, 
		\end{equation}
		where, from \eqref{conv rho eps}, \eqref{conv c eps} and \eqref{conv f eps},
		\begin{equation*}
			\sqrt{\vr_{\ve}} f_{\ve} c_{\ve} \rightharpoonup \overline{\sqrt{\vr} f c} \quad \mbox{in } L^2(0,T; L^{p_2}(\Omega)) \quad \mbox{with} \quad p_2= \frac{6\gamma}{4\gamma+3} > 1.  
		\end{equation*}
		\item From convergences \eqref{conv rho c eps}, \eqref{conv c eps} and the compact Sobolev embedding \eqref{SE2}, it follows that
		\begin{equation*}
			\vr_{\ve} c_{\ve}^2 \overset{*}{\rightharpoonup} \vr c^2 \quad \mbox{in } L^{\infty}(0,T; L^{p_4}(\Omega )) \quad \mbox{with} \quad p_4= \frac{3\gamma}{\gamma+3}>1.
		\end{equation*}
		Then, using additionally convergences \eqref{conv rho eps} and \eqref{conv rho c eps}, we deduce that 
		\begin{equation} \label{conv rho c2}
			\| \sqrt{\vr_{\ve}} (c_{\ve}-c)\|_{L^2((0,T) \times \Omega)}^2 =  \int_{0}^{T} \int_{\Omega} (\vr_\ve c_{\ve}^2 - 2 \vr_{\ve}c_{\ve} c + \vr_{\ve} c^2 ) \ \dx \dt \to 0 
		\end{equation}
		as $\ve \to 0$. Therefore, for any $\varphi \in L^{\infty}((0,T) \times \Omega)$ we have 
		\begin{equation*}
			\int_{0}^{T} \int_{\Omega} \left(\overline{\sqrt{\vr} f c} - \overline{\sqrt{\vr}f } \, c \right) \varphi \ \dx \dt = \lim_{\ve \to 0} 	\int_{0}^{T} \int_{\Omega} \sqrt{\vr_{\ve}} f_{\ve} (c_{\ve}-c) \varphi \ \dx \dt 
		\end{equation*}
		where, from \eqref{conv f eps}, 
		\begin{equation*}
			\begin{aligned}
				&\left| \int_{0}^{T} \int_{\Omega} \sqrt{\vr_{\ve}} f_{\ve} (c_{\ve}-c) \varphi \ \dx \dt  \right| \\
				&\quad \leq \| f_{\ve} \|_{L^2 ((0,T)\times \Omega)} \| \sqrt{\vr_{\ve}} (c_{\ve}-c)\|_{L^2((0,T)\times \Omega)} \| \varphi\|_{L^{\infty}((0,T)\times \Omega) } \\
				&\quad \leq C(\varphi) \| \sqrt{\vr_{\ve}} (c_{\ve}-c)\|_{L^2((0,T)\times \Omega)}  \to 0
			\end{aligned}
		\end{equation*}
		as $\ve \to 0$, leading to \eqref{I4}.
	\end{itemize}

Now, we point out that, once 
	\begin{equation*}
		\nabla_x c_{\varepsilon} \to \nabla_x c \quad \mbox{in } L^2((0,T) \times \Omega; \mathbbm{R}^3)
	\end{equation*} 
	is established, the correspondent strong convergence \eqref{strong conv c eps} can be deduced as a consequence of the generalized Poincar\'{e} inequalities reported in \eqref{PI} and in Lemma \ref{Generalized Poincare}; more precisely, there exists a positive constant $C$, depending only on $\gamma, M, \Omega$ and the constant in \eqref{E1}, such that 
    \begin{align*}
			&\int_{0}^{T} \| c_{\varepsilon} - c \|_{W^{1,2}(\Omega)}^2 \ \textup{d}x 
            \\
            &\quad  \leq C \int_{0}^{T} \| \nabla_x(c_{\varepsilon} - c) \|_{L^2(\Omega; \mathbbm{R}^3)}^2 \ \textup{d}x + C\int_{0}^{T} \int_{\Omega} \varrho_{\varepsilon}|c_{\varepsilon}- c | \ \textup{d}x \textup{d}t
            \\
            &\quad  \leq C \int_{0}^{T} \| \nabla_x(c_{\varepsilon} - c) \|_{L^2(\Omega; \mathbbm{R}^3)}^2 \ \textup{d}x 
            + C 
             \| \vr_{\ve} \|_{L^1((0,T)\times \Omega)}^\frac12
            \left( \int_{0}^{T} \int_{\Omega} \varrho_{\varepsilon}|c_{\varepsilon}- c |^2 \ \textup{d}x \textup{d}t\right)^\frac12.
		\end{align*}
	Then, due to \eqref{conv rho c2}, the right-hand side vanishes as $\varepsilon \to 0 $.

	 As a consequence of \eqref{strong conv c eps}, passing to suitable subsequences as the case may be, we infer in particular that 
	\begin{equation} \label{conv a.e. c eps}
		c_{\ve} \to c \quad \mbox{a.e. in } (0,T) \times \Omega.  
	\end{equation}
	Thus, from hypothesis \eqref{viscosity coefficients} and convergence \eqref{conv u eps}, we get \eqref{conv S eps}.
	
	\subsection{Strong convergence of the densities}
	
	So far, from estimate  \eqref{E18}, passing to a suitable subsequence as the case may be, we can only deduce that 
	\begin{equation} \label{conv p eps bis}
		p(\vr_{\ve}) \rightharpoonup \overline{p(\vr)} \quad \mbox{in } L^q((0,T) \times \Omega), 	
	\end{equation}
	with the exponent $q=q(\gamma)$ defined as in \eqref{E18}. Therefore, in order to conclude that $\overline{p(\vr)} = p(\vr)$, we need that 
	\begin{equation} \label{conv a.e. rho eps}
		\vr_{\ve} \to \vr \quad \mbox{a.e. } (0,T) \times \Omega.
	\end{equation}
	To this end, it is sufficient to show 
	\begin{equation} \label{strong conv rho eps}
		\vr_{\ve} \to \vr \quad \mbox{in } L^1((0,T) \times \Omega). 
	\end{equation}
	
	The proof of \eqref{strong conv rho eps} is definitely non-trivial but nowadays well-understood. The idea was first developed by Lions in \cite{Lio}, and later adapted by Feireisl in \cite{Fei1} to handle the case of non-constant viscosity coefficients; we report the key ingredients in Appendix \ref{SC density}.

    \subsection{Boundedness of the concentration. Part 1.} 
    \label{Bound Concentrations}
	
	Next, we prove some properties of the limit concentration $c$. By definition of the approximating potential $F_\ve$, we have
	\begin{equation*}
		F'_\ve(s)= F'_\ve(1-\ve) + F''_\ve(1-\ve)(s-(1-\ve)) \geq F''_\ve(1-\ve)(s-1), \quad \forall \, s >1.
	\end{equation*}
	Since $\varrho_\ve$ is non-negative, we infer that
	\begin{align*}
		\int_{\lbrace{c_\ve}>1 \rbrace} 
		\varrho_\ve |F_\ve''(1-\ve)|^2 (c_\ve-1)^2 \ \dx \dt
		&\leq 
		\int_{\lbrace{c_\ve}>1 \rbrace} 
		\varrho_\ve |F_\ve'(c_\ve)|^2 \ \dx \dt
		\\
		&\leq \int_{0}^{T} \int_{\Omega} 
		\varrho_\ve |F_\ve'(c_\ve)|^2 \ \dx \dt.
	\end{align*}
	In light of \eqref{E14}, we deduce that
	\begin{equation*}
		\int_{\lbrace{c_\ve}>1 \rbrace} 
		\varrho_\ve (c_\ve-1)^2 \ \dx \dt
		\leq 
		\frac{C}{|F_\ve''(1-\ve)|^2}.
	\end{equation*}
	Repeating the same argument for $s<-1$, it is easily seen that
	\begin{equation*}
		\int_{\lbrace{c_\ve}<-1 \rbrace} 
		\varrho_\ve (c_\ve+1)^2 \ \dx \dt
		\leq 
		\frac{C}{|F_\ve''(-1+\ve)|^2}.
	\end{equation*}
	Thus, by the symmetry of $F_\ve$, we have
	\begin{equation*}
		\int_{0}^{T} \int_{\Omega} \varrho_\ve (|c_\ve|-1)_{+}^2 \ \dx \dt
		\leq 
		\frac{C}{|F_\ve''(1-\ve)|^2}= \frac{C}{|F''(1-\ve)|^2}.
	\end{equation*}
	From the a.e. convergence of the concentrations \eqref{conv a.e. c eps} and densities \eqref{conv a.e. rho eps}, an application of the Fatou lemma entails that
	\begin{align*}
		\int_0^T \int_{\Omega} \varrho (|c|-1)_{+}^2 \ \dx \dt 
		&\leq \liminf_{\ve \to 0} \int_0^T \int_{\Omega}
		\varrho_\ve (|c_\ve|-1)_{+}^2 \ \dx \dt
		\\
		&\leq 
		\liminf_{\ve \to 0} 
		\frac{C}{|F''(1-\ve)|^2}=0.
	\end{align*}
	Thus, we conclude that
	\begin{equation}
		\Big| \Big\lbrace (t,x)\in (0,T) \times \Omega \ : \ \varrho(t,x)>0, \, |c(t,x)|>1 \Big\rbrace\Big|=0.
	\end{equation}
	In other words, 
	\begin{equation*}
		-1 \leq c \leq 1  \quad \mbox{a.e. in } \lbrace{ \varrho > 0 \rbrace} \subset (0,T) \times \Omega.
	\end{equation*}
	Furthermore, exploiting once again the pointwise convergences \eqref{conv a.e. c eps} and \eqref{conv a.e. rho eps}, it is not difficult to show that
    \begin{equation*}
		\varrho_\ve |F'_\ve(c_\ve)|^2 \to \varrho |\widetilde{F'}(c)|^2
		\quad \mbox{a.e. in } \lbrace{ \varrho > 0 \rbrace} \subset  (0,T) \times \Omega,
	\end{equation*}
	where
	\begin{equation*}
		\widetilde{F'}(s):= 
		\begin{cases}
			+ \infty \quad &\text{if} \quad s=1,\\
			F'(s) \quad &\text{if} \quad s \in (-1,1),\\
			-\infty \quad &\text{if} \quad s=-1.
		\end{cases}
	\end{equation*}
	Now, by the Fatou lemma and \eqref{E14}, we derive that
	\begin{align*}
		\int_{0}^T \int_{\Omega}
		\varrho |\widetilde{F'}(c)|^2 \mathbbm{1}_{\lbrace\varrho>0\rbrace} \ \dx \dt
		&\leq 
		\liminf_{\ve \to 0}
		\int_{0}^T \int_{\Omega}
		\varrho_\ve |F_\ve'(c_\ve)|^2 \mathbbm{1}_{\lbrace\varrho>0\rbrace} \ \dx \dt
		\\
		&\leq 
		\liminf_{\ve \to 0}
		\int_{0}^T \int_{\Omega}
		\varrho_\ve |F_\ve'(c_\ve)|^2 \ \dx \dt
		\leq C.
	\end{align*}
	As a consequence, we obtain that
	\begin{equation} \label{bound c}
		-1 < c < 1  \quad \mbox{a.e. in } \lbrace{ \varrho > 0 \rbrace} \subset (0,T)\times \Omega,
	\end{equation}
	which in turn implies that $\varrho F'(c)$ is well-defined on the set $\lbrace \varrho >0 \rbrace$.
	
	Owing to the above convergence results, we are now in the position to characterize the weak limit of $\varrho_\ve F'_\ve(c_\ve)$. To this end, we first recall that 
	\begin{equation}
		\label{strong-rho-0}
		\varrho_\ve \to 0 \quad \text{in} \quad L^1(\lbrace{\varrho=0\rbrace}). 
	\end{equation}
	Now, for any $\varphi \in C_c^\infty((0,T)\times \overline{\Omega})$, we have
	\begin{equation} \label{I5}
		\int_0^T \int_{\Omega} \varrho_\ve F'_\ve(c_\ve) \varphi \ \dx\dt =
		\int_{\lbrace{\varrho=0\rbrace} } \varrho_\ve F'_\ve(c_\ve) \varphi \ \dx\dt  
		+
		\int_{\lbrace{\varrho>0\rbrace} } \varrho_\ve F'_\ve(c_\ve) \varphi \ \dx\dt.  
	\end{equation}
	By \eqref{conv rho eps}, \eqref{E14} and \eqref{strong-rho-0}, we observe that
	\begin{align*}
		\left| \int_{\lbrace{\varrho=0\rbrace} } \varrho_\ve F'_\ve(c_\ve) \varphi \ \dx\dt   \right|
		&= \left| \int_{\lbrace{\varrho=0\rbrace} } \sqrt{\varrho_\ve} \sqrt{\varrho_\ve} F'_\ve(c_\ve) \varphi \ \dx\dt   \right|
		\\
		&\leq \| \sqrt{\varrho_\ve} \|_{L^2(\lbrace{\varrho=0\rbrace})} 
		\| \sqrt{\varrho_\ve} F'_\ve(c_\ve) \|_{L^2((0,T)\times \Omega} \| \varphi\|_{L^\infty( (0,T) \times \Omega}
		\\
		&\leq C \| \varrho_\ve \|_{L^1(\lbrace{\varrho=0\rbrace})}^\frac12 
		\to 0, \quad \text{for}\quad \ve \to 0.
	\end{align*}
	On the other hand, since
	\begin{equation*}
		\varrho_\ve F'_\ve(c_\ve) \to \varrho F'(c)
		\quad \mbox{a.e. in } \lbrace{ \varrho > 0 \rbrace} \subset (0,T) \times \Omega,
	\end{equation*}
	and, from estimates \eqref{E1} and \eqref{E14}, 
	\begin{equation}
		\|\varrho_\ve F'_\ve(c_\ve)\|_{L^2(0,T;L^{p_1}(\Omega))} \leq C \quad \mbox{with} \quad p_1= \frac{2\gamma}{\gamma+1}, 
	\end{equation}
	we have that 
	\begin{equation}
		\varrho_\ve F'_\ve(c_\ve) \mathbbm{1}_{\lbrace\varrho>0\rbrace} \rightharpoonup \varrho F'(c) \mathbbm{1}_{\lbrace\varrho>0\rbrace}
		\quad \mbox{in } L^2(0,T;L^\frac{2\gamma}{\gamma+1}(\Omega)).
	\end{equation}
	In conclusion, going back to \eqref{I5}, we deduce that
	\begin{equation}
		\int_0^T \int_{\Omega} \varrho_\ve F'_\ve(c_\ve) \varphi \ \dx\dt \to 
		\int_0^T \int_{\Omega} \varrho F'(c) \varphi \ \dx\dt, 
	\end{equation}
	where the function $\varrho F'(c)$ must be interpreted as
	$$
	\varrho F'(c)
	= 
	\begin{cases}
		\varrho F'(c) \quad &\text{if} \quad \varrho>0, \\
		0 \quad &\text{if} \quad \varrho=0.
	\end{cases}
	$$
	In particular, we get convergence \eqref{conv rho F' eps}.
	
	\subsection{Boundedness of the concentration. Part 2.}
	\label{Bound Concentrations_2}

We now take advantage of the uniform estimate \eqref{E14-bis} to infer the boundedness of the limit function $c$ on the set $\lbrace \varrho =0 \rbrace$. To this end, we consider $g(s)= \arcsin(s)$ in its domain $[-1,1]$ and we introduce its $C^1$-approximation 
$$
g_{\ve}(s):= \begin{cases}
			\displaystyle g(1-\ve) + g'(1-\ve)\left( s- (1-\ve) \right) &\mbox{for } 
			s > 1-\ve, 
            \\[8pt] 
            g(s) &\mbox{for } -1+\ve \leq s \leq 1-\ve, 
            \\[5pt]
			\displaystyle g(-1+\ve) + g'(-1+\ve)\left( s- (-1+\ve) \right)
			 &\mbox{for } s < -1+\ve,
		\end{cases}
$$
for any $\ve \in (0,1)$ sufficiently small. It is easily seen that $\sqrt{\theta} g'_\ve(s) = \sqrt{F''_\ve(s)}$, for any $s \in \mathbbm{R}$. Then, in light of \eqref{E14-bis}, we notice that
\begin{align*}
\int_0^T \int_{\Omega} |\Grad g_\ve(c_\ve)|^2 \ \dx\dt
= \int_0^T \int_{\Omega} |g'_\ve(c_\ve)|^2 |\Grad c_\ve|^2 \ \dx\dt
= \frac{1}{\theta} \int_0^T \int_{\Omega} F''_\ve(c) |\Grad c_\ve|^2 \ \dx\dt\leq C,
\end{align*}
namely 
\begin{equation} \label{G-1}
	\|  \Grad g_\ve(c_{\ve}) \|_{L^2(0,T; L^2(\Omega))} \leq C.
		\end{equation}
On the other hand, a simple calculation shows that 
$$
|g_\ve (s)| \leq \frac{2}{\theta} |F'_\ve(s)|, \quad \forall \, s \in \mathbbm{R}. 
$$
Hence, recalling \eqref{E14}, we immediately infer that 
\begin{align*}
\int_0^T \int_{\Omega} \varrho_\ve |g_\ve(c_\ve)| \ \dx\dt
&\leq \frac{2}{\theta} \int_0^T \int_{\Omega}  \varrho_\ve |F'_\ve(c_\ve)|
\\
&\leq \frac{2}{\theta} 
 \| \vr_{\ve} \|_{L^1((0,T)\times \Omega)}^\frac12
\left( \int_0^T \int_{\Omega}  \varrho_\ve |F'_\ve(c_\ve)|^2 \ \dx\dt\right)^\frac12
\leq C.
\end{align*}
Thus,  by exploiting \eqref{G-1} and the generalized Poincar\'{e} inequality in Lemma \ref{Generalized Poincare}, the latter entails that
$$
\| g_\ve (c_\ve) \|_{L^2(0,T; W^{1,2}(\Omega))} \leq C.
$$
Now, by a symmetry argument, we observe that 
$$
|g_\ve(s)|= g(1-\ve) + g'(1-\ve) \left( |s|-1 \right)  \geq g'(1-\ve) \left( |s|-1 \right), \quad \forall \, |s| >1.
$$
Therefore, we have
\begin{align*}
\int_0^T \int_\Omega \left( |c_\ve|-1 \right)_+^2 \ \dx\dt
\leq \frac{1}{|g'(1-\ve)|^2} \int_0^T \int_\Omega |g_\ve(c)|^2  \ \dx\dt
\leq \frac{C}{|g'(1-\ve)|^2}.
\end{align*}
From the a.e. convergence of the concentrations \eqref{conv a.e. c eps} and the Fatou lemma, we deduce that 
$$
\int_0^T \int_\Omega \left( |c|-1 \right)_+^2 \ \dx\dt
\leq \liminf_{\ve \to 0} \int_0^T \int_\Omega \left( |c_\ve|-1 \right)_+^2 \ \dx\dt
\leq \liminf_{\ve \to 0} \frac{C}{|g'(1-\ve)|^2}=0,
$$
which implies that 
$$
-1 \leq c \leq 1  \quad \mbox{a.e. in } (0,T) \times \Omega.
$$

   \subsection{Passage to the limit}
	
	Thanks to the necessary convergences listed in Lemma \ref{Lemma Conv}, we are now ready to let $\ve \to 0$ in identities  \eqref{WF2}--\eqref{WF5}, obtaining that the weak formulations \eqref{WF2fin}--\eqref{WF5fin} hold true for the limit functions $(\vr, \vu, c)$; we point out, in particular, that the validity of the renormalized continuity equation \eqref{WF2} was achieved in Section \ref{Ren}, when proving the strong convergence of the densities. 
	
	In order to conclude the proof Theorem \ref{Main Result}, it remains to show that the energy inequality \eqref{energy inequality final} holds true for a.e. $\tau \in (0,T)$. To this end, we first observe that, due to the convexity of $F_{\ve}=F_{\ve}(s)$, 
	\begin{equation*}
		F_{\ve}(s_0) \geq F_{\ve}(s) + F_{\ve}'(s)(s_0-s) \quad \mbox{for any } s, s_0 \in \mathbbm{R}.
	\end{equation*}
	Hence, choosing $s=c_{\ve}$ and $s_0=0$, we obtain that $F_{\ve}(c_{\ve}) \leq c_{\ve} F_{\ve}'(c_\ve) $, which implies, in view of estimates \eqref{E1}, \eqref{E6} and \eqref{E14} that 
	\begin{equation*}
		\| \vr_{\ve} F_{\ve} (c_{\ve})\|_{L^2(0,T; L^{p_2}(\Omega))} \leq C  \quad \mbox{with} \quad p_2= \frac{6\gamma}{4\gamma +3}. 
	\end{equation*}
	Proceeding similarly as in Section \ref{Bound Concentrations}, from the boundedness of the concentrations \eqref{bound c} and the pointwise convergences \eqref{conv a.e. c eps} and \eqref{conv a.e. rho eps}, we have  
	\begin{equation*}
		\vr_{\ve} F(c_{\ve}) \to \vr F(c) \quad \mbox{a.e. on } \{ \vr>0 \} \subset (0,T) \times \Omega.
	\end{equation*}
	On the other hand, from the strong convergence \eqref{strong-rho-0} of the densities on the set $\{\vr=0\}$, for any $\varphi \in C_c^{\infty}((0,T) \times \Omega)$ we have 
	\begin{equation*}
		\begin{aligned}
			\left| \int_{\{\vr=0\}} \vr_{\ve} F_{\ve}(c_{\ve}) \varphi \ \dx \dt  \right| &\leq C(\varphi) \| \sqrt{\vr_{\ve}}\|_{L^3(\{\vr=0\})} \| \sqrt{\vr_{\ve}} F_{\ve}'(c_{\ve}) \|_{L^2((0,T) \times \Omega)} \| c_{\ve} \|_{L^{\infty}(0,T; L^6(\Omega))} \\
			&\leq C \| \sqrt{\vr_{\ve}}\|_{L^1(\{\vr=0\})}^{\frac12} \to 0 \quad \mbox{as} \quad \ve \to 0.  
		\end{aligned} 
	\end{equation*}
	We can therefore conclude that 
	\begin{equation*}
		\vr_{\ve} F(c_{\ve}) \rightharpoonup \vr F(c) \mathbbm{1}_{\{ \vr>0\}} \quad \mbox{in } L^2(0,T; L^{\frac{6\gamma}{4\gamma +3}}(\Omega)). 
	\end{equation*}
	Furthermore, in light of the strong convergences \eqref{strong conv c eps}, \eqref{strong conv rho eps}, it is not difficult to prove that 
	\begin{equation*}
		\vr_{\ve} c_{\ve}^2 \to \vr c^2  \quad \mbox{in } L^1((0,T) \times \Omega );
	\end{equation*}
	indeed, using additionally estimates \eqref{E1}, \eqref{E6}, the Sobolev embedding \eqref{SE1} and interpolation, we have, for some $\omega \in (0,1)$,
	\begin{equation*}
		\begin{aligned}
			\| \vr_{\ve} c_{\ve}^2 - \vr c^2 \|_{L^1((0,T) \times \Omega )} &\leq  \int_{0}^{T} \int_{\Omega} \vr_{\ve} |c_{\ve}+ c| \, |c_{\ve} - c| \ \dx \dt + \int_{0}^{T} \int_{\Omega}  |\vr_{\ve}- \vr| \, c^2 \ \dx \dt \\
			 & \leq 2 \| \vr_{\ve}\|_{L^{\infty}(0,T; L^\frac32(\Omega))}\| c_{\ve}\|_{L^2(0,T; L^6(\Omega))} \| c_{\ve} - c \|_{L^2(0,T; L^6(\Omega))} \\
			 &+ \| c\|^2_{L^6((0,T) \times \Omega)} \| \vr_{\ve}- \vr \|_{L^{\frac32}((0,T) \times \Omega )} \\
			 &\leq C\left( \| c_{\ve} - c \|_{L^2(0,T; W^{1,2}(\Omega))} + \| \vr_{\ve}- \vr \|_{L^1((0,T) \times \Omega)}^{\omega} \right) \to 0 \quad \mbox{as} \quad \ve \to 0.   
		\end{aligned}
	\end{equation*}
	Hence, exploiting the lower semicontinuity of the norms with respect to the weak or weak-$*$ convergence, we are allowed to pass to the limit in the energy inequality \eqref{energy inequality}; notice that it is first necessary to integrate $E_{\ve}= E_{\ve}(t)$ over $(\tau, \tau +h)$  to deduce 
	\begin{equation} \label{e1}
		\int_{\tau}^{\tau+h} E(t) \ \dt  \leq \liminf_{\ve \to 0} 	\int_{\tau}^{\tau+h} E_{\ve}(t)  \ \dt 
	\end{equation}
	for any $h>0$, where 
	\begin{equation*}
		E(t):= \int_{\Omega} \left( \frac12 \vr |\vu|^2 + \frac{1}{\gamma-1}\vr^{\gamma} + \vr F(c) \mathbbm{1}_{\{ \vr>0\}}  - \frac{\theta_0}{2}\vr c^2   + \frac12 |\Grad c|^2\right)(t, \cdot) \ \dx. 
	\end{equation*}
	Next, it is enough to divide the left-hand side of \eqref{e1} by $h$ and let $h\to 0$ to recover that the energy inequality holds for a.e. $\tau \in (0,T)$.

    \appendix

    \section{Strong convergence of the densities} \label{SC density}
	
	In this section we report the main steps in order to prove the strong convergence \eqref{strong conv rho eps} of the densities. 
	\subsection{Step 1}
		Introducing the ``truncation" function, defined for any $k \in \mathbb{N}$ as 
		\begin{equation*}
			T_k(\vr):= \min\{ \vr, k\}, 
		\end{equation*}
		the first step consists in showing that 
		\begin{equation} \label{Step 1}
			\begin{aligned}
				\lim_{\ve \to 0} &\int_{0}^{T} \psi \int_{\Omega} \phi \, \Big( p(\vr_{\ve}) \, \xi \, T_k (\vr_{\ve}) - \vS (c_{\ve}, \Grad \vu_{\ve}) : \Grad \Delta_x^{-1} \Grad \big[ \xi \, T_k(\vr_{\ve})\big] \Big) \dx \dt \\
				=& \ \int_{0}^{T} \psi \int_{\Omega} \phi \, \Big( \overline{p(\vr)} \, \xi \, \overline{T_k (\vr)} - \vS (c, \Grad \vu) : \Grad \Delta_x^{-1} \Grad \big[ \xi \, \overline{T_k(\vr)}\big] \Big) \dx \dt, 
			\end{aligned}
		\end{equation}
		for any $\psi \in C_c^{\infty}(0,T)$ and any $\phi, \xi \in C_c^{\infty}(\Omega)$. Here, the symbol $\Delta_x^{-1}$ denotes the inverse of the Laplace operator considered on the whole space $\mathbbm{R}^3$, while 
		\begin{equation} \label{conv T_k rho}
			T_k(\vr_{\ve}) \to \overline{T_k(\vr)} \quad \mbox{in } C_{\rm weak} ([0,T]; L^p(\Omega)) \ \mbox{for any finite } p>1. 
		\end{equation}
		Due to the convergences \eqref{conv rho eps}--\eqref{conv S eps} along with \eqref{conv p eps bis}, we can pass to the limit in \eqref{WF3}, obtaining that the integral identity 
		\begin{equation} \label{WF3bis}
			\begin{aligned}
				\left[\int_{\Omega} (\vr \vu \cdot \bm{\varphi})(t, \cdot) \ \dx\right]_{t=0}^{t=\tau} &=\int_{0}^{\tau} \int_{\Omega} \big[ \vr \vu \cdot \DerTime \bm{\varphi} +(\vr \vu \otimes \vu) : \Grad \bm{\varphi}  \big] \, \dx \dt \\
				&+ \int_{0}^{\tau} \int_{\Omega} \Big[\overline{p(\vr)} \mathbb{I}- \vS(c, \Grad \vu)\Big] : \Grad \bm{\varphi} \ \dx \dt \\
				&+ \int_{0}^{\tau} \int_{\Omega} \left( (\Grad c \otimes \Grad c ) : \Grad \bm{\varphi} - \frac{|\Grad c|^2}{2} \Div \bm{\varphi}\right) \dx \dt
			\end{aligned}
		\end{equation}
		holds for any for any $\tau \in [0,T]$ and any $\bm{\varphi} \in C^\infty ([0,T) \times \Omega; \mathbbm{R}^3 )$. Hence, we can take the quantities 
		\begin{align*}
			\bm{\varphi}_{\ve} (t,x) &:= \psi(t) \phi(x) \Grad \Delta_x^{-1} \left[ \xi(x) T_k(\vr_{\ve})(t,x)\right], \\ 
			\bm{\varphi}(t,x) &:= \psi(t) \phi(x) \Grad \Delta_x^{-1} \left[ \xi(x) \overline{T_k(\vr)}(t,x)\right], 
		\end{align*}
		with $\psi \in C^{\infty}_c(0,T)$ and $\phi, \, \xi \in C_c^{\infty}(\Omega)$ fixed, as test functions in \eqref{WF3} and \eqref{WF3bis}, respectively. Identity \eqref{Step 1} will then follow by letting $\ve \to 0$, observing that $T_k(\vr_{\ve})$ satisfies the renormalized continuity equation in the sense of distributions and exploiting the boundedness of the singular operator
		\begin{equation} \label{boundedness singular operator}
			\Grad \Delta_x^{-1} : L^p(\Omega) \to W^{1,p}(\Omega) \quad \mbox{for any } p>1.
		\end{equation}
		The main tool needed in this passage is a slightly modified version of the ``Div--Curl Lemma" contained in \cite[Corollary 6.1, (ii)]{Fei} to deduce that 
		\begin{equation*}
			\vv_{\ve} (t)  \rightharpoonup \vv(t) \quad \mbox{in } L^r(\Omega; 	\mathbbm{R}^3) \quad \mbox{for any } 1\leq r < \frac{2\gamma}{\gamma +1}, 
		\end{equation*}
		where the quantities 
		\begin{align*}
			\vv_{\ve} (t) &:= \Grad \Delta_x^{-1} \Grad \left[ \xi \,  T_k(\vr_{\ve} ) (t, \cdot )\right] \phi \, (\vr_{\ve} \vu_{\ve})(t, \cdot)  - \xi \, T_k(\vr_{\ve})(t, \cdot) \Grad \Delta_x^{-1} \Div \left[ \phi \, (\vr_{\ve} \vu_{\ve})(t, \cdot) \right], \\
			\vv(t) &:= \Grad \Delta_x^{-1} \Grad \left[ \xi \,  \overline{T_k(\vr)}(t, \cdot) \right] \phi \, (\vr \vu)(t, \cdot) - \xi \, \overline{T_k(\vr)}(t, \cdot) \Grad \Delta_x^{-1} \Div \left[ \phi \, (\vr \vu)(t, \cdot) \right]
		\end{align*}
		are defined for any $t \in [0,T]$ -- we point out that all the functions involved are continuous in time.
		
		\subsection{Step 2}
		
		For a fixed $k \in  \mathbb{N}$, the next step is to show that 
		\begin{equation} \label{Step 2}
			\overline{p(\vr) \, T_k(\vr)} - \overline{p(\vr)} \, \overline{T_k(\vr)}= \left(\lambda(c) + \frac43 \mu(c)\right) \left( \overline{T_k(\vr) \Div \vu} - \overline{T_k(\vr)} \Div \vu \right). 
		\end{equation}
		
		We observe that, for any vector $\vv$,
		\begin{equation*}
			\Grad \Delta_x^{-1} \Grad : [\Grad \vv] = \sum_{j,k=1}^{3}  \mathcal{F}^{-1} \left( \frac{\xi_j \xi_k}{|\xi|^2} \mathcal{F}\left( \partial_{x_k} v_j\right)\right) = \sum_{j=1}^{3} \mathcal{F}^{-1} \left( i\xi_j \mathcal{F}\left(  v_j\right)\right) = \Div \vv, 
		\end{equation*}
		and therefore we can write 
		\begin{equation*}
			\begin{aligned}
				&\int_{0}^{T} \psi \int_{\Omega} \phi \, \vS (c_{\ve}, \Grad \vu_{\ve}) : \Grad \Delta_x^{-1} \Grad \big[ \xi \, T_k(\vr_{\ve})\big]  \dx \dt \\
				=& \,2 \int_{0}^{T} \psi \int_{\Omega} \xi \, \Big(  \Grad \Delta_x^{-1} \Grad :  \big[ \phi \,\mu(c_{\ve}) \Grad \vu_{\ve}\big]  - \phi \, \mu(c_{\ve})   \Grad \Delta_x^{-1} \Grad : [\Grad \vu_{\ve}]  \Big) \, T_k(\vr_{\ve}) \  \dx \dt \\
				&+ \int_{0}^{T} \psi \int_{\Omega} \xi \, \phi  \left(\lambda(c_{\ve}) + \frac43 \mu(c_{\ve})\right) T_k(\vr_{\ve}) \,  \Div \vu_{\ve}   \  \dx \dt.
			\end{aligned}
		\end{equation*}
		By interpolation we have 
		\begin{equation*}
			L^{\infty}(0,T; W^{1,2}(\Omega)) \cap L^2(0,T; W^{1,2p_4}(\Omega)) \hookrightarrow L^{2p_2}(0,T; W^{1,2p_2}(\Omega)) \quad \mbox{with} \quad p_2= \frac{6\gamma}{4\gamma+3}> \frac65, 
		\end{equation*}
		and therefore, proceeding as for \eqref{E16}, we can apply a variant of the ``Commutator Lemma", see \cite[Lemma 4.2]{Fei1}, to deduce that 
		\begin{equation*}
			\left\| \Grad \Delta_x^{-1} \Grad :  \big[ \phi \,\mu(c_{\ve}) \Grad \vu_{\ve}\big]  - \phi \, \mu(c_{\ve})   \Grad \Delta_x^{-1} \Grad : [\Grad \vu_{\ve}] \right\|_{L^q(0,T; W^{\alpha, q}(\Omega))} \leq C
		\end{equation*}
		for some $q>1$ and $0< \alpha <1$. Hence, from \eqref{conv u eps},  \eqref{conv a.e. c eps} and \eqref{conv T_k rho} with $p$ large enough so that $\frac1p+\frac1q<1$, we deduce 
		\begin{equation*}
			\begin{aligned}
				\lim_{\ve \to 0} &\int_{0}^{T} \psi \int_{\Omega} \xi \, \Big(  \Grad \Delta_x^{-1} \Grad :  \big[ \phi \,\mu(c_{\ve}) \Grad \vu_{\ve}\big]  - \phi \, \mu(c_{\ve})   \Grad \Delta_x^{-1} \Grad : [\Grad \vu_{\ve}]  \Big) \, T_k(\vr_{\ve}) \  \dx \dt \\
				=& \int_{0}^{T} \psi \int_{\Omega} \xi \, \Big(  \Grad \Delta_x^{-1} \Grad :  \big[ \phi \,\mu(c) \Grad \vu\big]  - \phi \, \mu(c)   \Grad \Delta_x^{-1} \Grad : [\Grad \vu]  \Big) \, \overline{T_k(\vr)} \  \dx \dt. 
			\end{aligned}
		\end{equation*}
		Going back to \eqref{Step 1}, we finally get \eqref{Step 2}.
		
		\subsection{Step 3} \label{Ren}
		
		We now need to show that the limits $(\vr, \vu)$ of \eqref{conv rho eps}, \eqref{conv u eps} constitute a renormalized solution, in the sense that the integral identity \eqref{WF2fin} holds for any  given functions $b$ and $B$ satisfying \eqref{b and B}. In view of the regularizing procedure by DiPerna and Lions \cite{DiPLio}, see also 
		\cite[Corollary 4.1]{Fei}, the latter is satisfied if the density $\vr$ is at least square-integrable, which, due to estimate \eqref{bound rho}, holds true for $\gamma \geq \frac95$. To handle the case $\frac32 < \gamma < \frac95$, as consequence of \cite[Proposition 6.3]{Fei} it is enough to show that the \textit{oscillation defect measure} is finite:
		\begin{equation} \label{oscillation defect measure}
			\sup_{k \in \mathbb{N}} \left(\limsup_{\ve \to 0} \int_{0}^{T} \int_{\Omega} | T_k(\vr_{\ve}) - T_k(\vr)|^{\gamma+1} \dx \dt \right) <\infty. 
		\end{equation}
		From the convexity of $p=p(\vr)$, identity \eqref{Step 2}, Young inequality and hypothesis \eqref{viscosity coefficients}, we have
		\begin{equation*}
			\begin{aligned}
			&\limsup_{\ve \to 0} \int_{0}^{T} \int_{\Omega} | T_k(\vr_{\ve}) - T_k(\vr)|^{\gamma+1} \dx \dt \\
			& \leq \lim_{\ve \to 0} \int_{0}^{T} \int_{\Omega} \left( p(\vr_{\ve}) T_k(\vr_{\ve}) - \overline{p(\vr)} \, \overline{T_k(\vr)}\right) \dx \dt \\
			&= \lim_{\ve \to 0} \int_{0}^{T} \int_{\Omega} \left(\lambda(c) + \frac43 \mu(c)\right) T_k(\vr_{\ve})  \left(  \Div \vu_{\ve} - \Div \vu \right) \dx\dt \\
			&\leq  \frac{1}{\gamma+1}   \limsup_{\ve \to 0} \int_{0}^{T} \int_{\Omega} | T_k(\vr_{\ve}) - T_k(\vr)|^{\gamma+1} \dx \dt +\frac{\gamma}{\gamma+1}\lim_{\ve \to 0} \int_{0}^{T} \int_{\Omega} | \Div(\vu_{\ve}-\vu)|^{1+\frac{1}{\gamma}} \ \dx \dt,  
			\end{aligned}
		\end{equation*}
		and hence, \eqref{oscillation defect measure} follows from \eqref{conv u eps}. 
		
		\subsection{Step 4}
		
		Finally, we choose
		\begin{equation*}
			b_k(\vr) \equiv T_k(\vr):= \min\{ \vr, k\}, \quad B_k(\vr) \equiv L_k(\vr):= \int_{1}^{\vr} \frac{ \min\{ z, k\}}{z^2} \ \textup{d}z, 
		\end{equation*}
		and $\varphi \equiv 1$ in both \eqref{WF2fin} and \eqref{WF2}. By using \eqref{Step 2}, we find that the identity
		\begin{equation}
			\begin{aligned}
				\int_{\Omega} \left( \overline{\vr L_k(\vr)} - \vr L_k(\vr) \right)(\tau, \cdot) \ \dx &+ \int_{0}^{\tau} \int_{\Omega} \left(\lambda(c) + \frac43 \mu(c)\right)^{-1} \left(	\overline{p(\vr) \, T_k(\vr)} - \overline{p(\vr)} \, \overline{T_k(\vr)}\right)  \dx \dt \\
				&+\int_{0}^{\tau} \int_{\Omega} \left( \overline{T_k(\vr)} - T_k(\vr) \right)\Div \vu \  \dx \dt   =0  
			\end{aligned}
		\end{equation}
		holds for any $\tau \in [0,T]$. The last integral vanishes for $k\to +\infty$; indeed, from \eqref{oscillation defect measure}, H\"older inequality and interpolation, we have that
	\begin{equation*}
		\begin{aligned}
			&\left| \int_{0}^{\tau} \int_{\Omega} \left( \overline{T_k(\vr)} - T_k(\vr) \right)\Div \vu \  \dx \dt  \right| \\
			&\leq \|  \overline{T_k(\vr)} - T_k(\vr) \|_{L^2((0,T) \times \Omega)} \| \Div \vu \|_{L^2((0,T) \times \Omega)} \\
			&\leq \|  \overline{T_k(\vr)} - T_k(\vr) \|_{L^1((0,T) \times \Omega)}^{\omega} \|  \overline{T_k(\vr)} - T_k(\vr) \|_{L^{\gamma+1}((0,T) \times \Omega)}^{1-\omega} \| \Div \vu \|_{L^2((0,T) \times \Omega)},
		\end{aligned}
	\end{equation*}
	for some $\omega \in (0,1)$; from the weak lower semi-continuity of the norm, we deduce 
	\begin{equation*}
		\begin{aligned}
			\|  &\overline{T_k(\vr)} - T_k(\vr) \|_{L^1((0,T) \times \Omega)} \\
			& \leq \liminf_{\ve \to 0} \|  T_k(\vr_{\ve}) - \vr_{\ve}\|_{L^1((0,T) \times \Omega)} + \|  \vr  - T_k(\vr) \|_{L^1((0,T) \times \Omega)} \\
			&\leq \sup_{\ve>0} \int_{\{ \vr_{\ve} \geq k\}} \vr_{\ve} \ \dx \dt  + \int_{\{ \vr \geq k\}} \vr \ \dx \dt  \\
			&\leq k \sup_{\ve>0} \left( \int_{\left\{ \frac{\vr_{\ve}}{k} \geq 1 \right\}}\left(\frac{ \vr_{\ve}}{k}\right)^\gamma \ \dx \dt  + \int_{\left\{ \frac{\vr}{k} \geq 1 \right\}}\left(\frac{ \vr}{k}\right)^\gamma \ \dx \dt \right) \leq 2 \, k^{1-\gamma} \sup_{\ve>0} \| \vr_{\ve}\|_{L^{\gamma}((0,T) \times \Omega)}^{\gamma},
		\end{aligned}
	\end{equation*}
	where the last term on the right-hand side vanishes for $k \to +\infty$.
	
	Moreover, from hypothesis \eqref{viscosity coefficients}, there exist constants $\underline{V}, \, \overline{V}$ such that 
	\begin{equation*}
		0 <\underline{V} \leq \left(\lambda(c) + \frac43 \mu(c)\right)^{-1} \leq \overline{V}, 
	\end{equation*}
	while, from the fact that the function $p=p(\vr)$ is non-decreasing in $[0,\infty)$, we have that 
	\begin{equation*}
		\int_{0}^{\tau} \int_{\Omega} \left(	\overline{p(\vr) \, T_k(\vr)} - \overline{p(\vr)} \, \overline{T_k(\vr)}\right)  \dx \dt \geq 0. 
	\end{equation*}
	Thus, we obtain, for any $\tau \in [0,T]$,
	\begin{equation*}
		\int_{\Omega}  \left( \overline{\vr \ln(\vr)} - \vr \ln(\vr) \right)(\tau, \cdot) \ \dx \leq C \lim_{k \to \infty}\int_{\Omega} \left( \overline{\vr L_k(\vr)} - \vr L_k(\vr) \right)(\tau, \cdot) \ \dx \leq 0; 
	\end{equation*}
	on the other side, $\vr \ln(\vr) \leq \overline{\vr \ln(\vr)}$ due to the convexity of the function $\vr \mapsto \vr \ln (\vr)$ and therefore we get that 
	\begin{equation*}
		\overline{\vr \ln(\vr)} = \vr \ln (\vr ) \quad \mbox{a.e. on } (0,T) \times \Omega,
	\end{equation*}
	leading to \eqref{strong conv rho eps}.

\bigskip

\noindent
\textbf{Acknowledgments.} 
The authors wish to thank the anonymous referees for the careful reading and helpful suggestions.
The authors are supported by the MUR grant Dipartimento di Eccellenza 2023-2027 of Dipartimento di Matematica, Politecnico di Milano. 
D. Basari\'{c} is supported by the INdAM-GNAMPA project ``EDP e applicazioni: dinamica dei fluidi e teoria spettrale", CUP E53C25002010001, and by the PRIN project 2022 ``Partial differential equations and related geometric-functional inequalities'', financially supported by the EU, in the framework of the ``Next Generation EU initiative''.
A. Giorgini is supported by the INdAM-GNAMPA project ``Analisi di modelli di Cahn-Hilliard per la separazione di fase'', CUP E53C25002010001.

\bigskip

\noindent
\textbf{Disclosure statement.} The authors report there are no competing interests to declare.

\medskip

\noindent
\textbf{Data availability statement.} No further data is used in this manuscript.

\end{document}